\RequirePackage{fix-cm}
\documentclass[smallextended]{svjour3}       
\smartqed  
\usepackage{graphicx}
\usepackage{amssymb}
\usepackage{mathrsfs}
\usepackage{amsmath}
\usepackage{float}
\usepackage[breaklinks,bookmarks=false]{hyperref}
\usepackage{enumitem}

\hypersetup{colorlinks, linkcolor=blue, citecolor=blue,
  urlcolor=blue, pagecolor=blue}

\floatstyle{ruled}
\newfloat{algorithm}{tbp}{loa}
\providecommand{\algorithmname}{Algorithm}
\floatname{algorithm}{\protect\algorithmname}

\usepackage{algorithmic}

\renewcommand{\S}{Section~}

\usepackage{thmtools}

\newcommand\mA{\mathbb{A}}

\newcommand\mP{\mathbb{P}}

\newcommand\krl{{\mathscr L}}

\newcommand\krr{{\mathscr R}}

\newcommand\krt{{\mathscr T}}

\newcommand\Th{{\krt_h}}

\newcommand\Thm{{\krt_{h,m}}}

\newcommand\bx{{\bf x}}
\newcommand\by{{\bf y}}

\newcommand\bkb{\mbox{\boldmath$b$\unboldmath}}
\newcommand\bkc{\mbox{\boldmath$c$\unboldmath}}

\newcommand\bkn{\mbox{\boldmath$n$\unboldmath}}

\newcommand\bkp{\mbox{\boldmath$p$\unboldmath}}
\newcommand\bkr{\mbox{\boldmath$r$\unboldmath}}
\newcommand\bkq{\mbox{\boldmath$q$\unboldmath}}
\newcommand\bks{\mbox{\boldmath$s$\unboldmath}}

\newcommand\bkx{\mbox{\boldmath$x$\unboldmath}}
\newcommand\bky{\mbox{\boldmath$y$\unboldmath}}

\newcommand\R{{\mathbb R}}
\newcommand\IR{{\mathbb R}}

\newcommand\dS{\,{\rm d}S}

\newcommand\Om{\Omega}

\newcommand\gom{{\Gamma}}

\newcommand\norm[2]{{\left\|{#1}\right\|_{#2}^{}} }

\newcommand\ah {{a_h}}

\newcommand{\DoF} {{\mathrm{DoF}}}

\newcommand\ve {\varepsilon}

\newcommand{\vp} {{\varphi}}

\newcommand\errSk{{ e_{\mathrm{S},k}}}
\newcommand\errAk{{ e_{\mathrm{A},k}}}
\newcommand\errDSk{{ e_{\mathrm{S},k}^*}}
\newcommand\errDAk{{ e_{\mathrm{A},k}^*}}

\newcommand\errSmk{{ e_{\mathrm{S},m,k}}}
\newcommand\errAmk{{ e_{\mathrm{A},m,k}}}
\newcommand\errDSmk{{ e_{\mathrm{S},m,k}^*}}
\newcommand\errDAmk{{ e_{\mathrm{A},m,k}^*}}

\newcommand\estAk {\zeta_{\mathrm{A},k}}
\newcommand\estSk {\zeta_{\mathrm{S},k}}
\newcommand\estDAk {\zeta_{\mathrm{A},k}^*}
\newcommand\estDSk {\zeta_{\mathrm{S},k}^*}

\newcommand\estAmk {\zeta_{\mathrm{A},m,k}}
\newcommand\estSmk {\zeta_{\mathrm{S},m,k}}
\newcommand\estDAmk {\zeta_{\mathrm{A},m,k}^*}
\newcommand\estDSmk {\zeta_{\mathrm{S},m,k}^*}

\newcommand\ESTAmk {\sigma_{\mathrm{A},m,k}}
\newcommand\ESTDAmk {\sigma_{\mathrm{A},m,k}^*}

\newcommand\etaAk {\eta_{\mathrm{A},k}}

\newcommand\etaSk {\eta_{\mathrm{S},k}}

\newcommand\etaSmk {\eta_{\mathrm{S},m,k}}
\newcommand\etaAmk {\eta_{\mathrm{A},m,k}}
\newcommand\etamk {\eta_{m,k}}
\newcommand\etaDSmk {\eta_{\mathrm{S},m,k}^*}
\newcommand\etaDAmk {\eta_{\mathrm{A},m,k}^*}
\newcommand\etaDmk {\eta_{m,k}^*}

\newcommand\J {J}
\newcommand\V {V}
\newcommand\Vh {V_h}
\newcommand\Wh {W_h}
\newcommand\Vhp {V_h^+}
\newcommand\Nh {{N_h}}

\newcommand\Vhm {V_{h,m}}
\newcommand\Vhz {V_{h,0}}
\newcommand\VhmP {V_{h,{m+1}}}

\renewcommand\Thm {\krt_{h,m}}
\newcommand\Thz {\krt_{h,0}}
\newcommand\ThmP {\krt_{h,{m+1}}}

\newcommand\uhm {u_{h,m}}
\newcommand\uhmk {u_{h,m}^k}
\newcommand\huhm {u_{h,m}^+}

\newcommand\zhm {z_{h,m}}
\newcommand\zhmk {z_{h,m}^k}
\newcommand\hzhm {z_{h,m}^+}

\newcommand\zhp {{z_h^+}}

\newcommand\TOL {{{\omega}}}
\newcommand\TOLA {{\omega_\mathrm{{A}}}}

\def \u {{u}} 
\def \z {{z}} 
\def \uh {{u_h}}
\def \zh {{z_h}}
\def \uhk {{u_{h}^{k}}}

\def \tuh {{\tilde{u}_h}}

\def \huh {{{u}_h^+}}

\def \hzh {{{z}_h^+}}

\def \zhk {{z_{h}^{k}}}

\providecommand{\res}[2]{ r_h(#1)(#2) }
\providecommand{\resD}[2]{ r_h^*(#1)(#2) }

\providecommand{\intX}[4]{\int_{#1}^{#2}#3\,\mathrm{d}#4}

\newcommand\restr[2]{{
  \left.\kern-\nulldelimiterspace 
  #1 
  \vphantom{\big|} 
  \right|_{#2} 
  }}

\def \etaDAk {\eta_{\mathrm{A},k}^{*}}

\def \etaDSk {\eta_{\mathrm{S},k}^{*}}

\newcommand\T{{\mkern-1.5mu\mathsf{T}}}
\newcommand\x{x}
\newcommand\y{y}
\renewcommand\c{c}
\renewcommand\b{b}

\renewcommand\bx{\bkx}
\renewcommand\by{\bky}
\newcommand\bxk{\bx_k}
\newcommand\byk{\by_k}

\newcommand\bykT{\by_k^{\T}}
\newcommand\bxknu{\bx_{k+\nu}}
\newcommand\byknu{\by_{k+\nu}}

\newcommand\byknuT{\by_{k+\nu}^{\T}}

\newcommand\bc{\bkc}
\newcommand\bb{\bkb}
\newcommand\bcT{\bkc^{\T}}

\newcommand\bp{\bkp}
\newcommand\bq{\bkq}
\newcommand\br{\bkr}
\newcommand\bs{\bks}
\newcommand\brk{\bkr_k}
\newcommand\bsk{\bks_k}

\newcommand\bskT{\bks_k^{\T}}

\newcommand\brknu{\bkr_{k+\nu}}

\newcommand\bsknuT{\bks_{k+\nu}^{\T}}

\newcommand\bxi{{\xi}}

\newcommand\bxiP{{\xi}^{\mathrm{p}}}
\newcommand\bxiD{{\xi}^{\mathrm{d}}}
\newcommand\bxiB{{\xi}^{\mathrm{B}}}
\newcommand\bxiBk{{\xi}^{\mathrm{B}}_k}
\newcommand\bxiBknu{{\xi}^{\mathrm{B}}_{k+\nu}}

\newcommand\PPa{{(P1)}}
\newcommand\PPb{{(P2)}}
\newcommand\PPc{{(P3)}}

\newcommand\EEa{{(E1)}}
\newcommand\EEb{{(E2)}}
\newcommand\EEc{{(E3)}}
\newcommand\cA{{c_{\mathrm{A}}}}

\usepackage[usenames]{color}
\newcommand{\Vit}[1]{{\color{red}{VD: #1} }}
\newcommand{\vit}[1]{{\color{red}{#1}}}

\newcommand{\Smaz}[1]{{}}

\usepackage{color}

\begin{document}

\title{On efficient numerical solution of linear algebraic systems arising
  in goal-oriented error estimates\thanks{This work was supported by grant No. 17-04150J
  of the Czech Science Foundation.}
}

\titlerunning{Solution of linear algebraic systems 
  in goal-oriented error estimates}        

\author{V{\'\i}t Dolej{\v s}{\'\i} \and Petr Tich{\'y}
}


\institute{V. Dolej{\v s}{\'\i}, P. Tich{\'y}\at
  Charles University Prague, \\
  Faculty of Mathematics and Physics,
  Sokolovsk\'a 83, 186 75 Praha, \\
  Czech Republic
  \email{dolejsi@karlin.mff.cuni.cz}          
  \email{ptichy@karlin.mff.cuni.cz}          
}

\date{Received: date / Accepted: date}

\maketitle

\begin{abstract}
  We deal with the numerical solution of linear partial differential equations (PDEs)
  with focus on the goal-oriented error estimates including algebraic errors
  arising by an inaccurate solution of the corresponding algebraic systems.
  The goal-oriented error estimates require the solution of the primal as well as
  dual algebraic systems. We solve both systems simultaneously using the
  bi-conjugate gradient method which allows to control the algebraic errors of both
  systems. We develop a stopping criterion which is cheap to evaluate and guarantees
  that the estimation of the algebraic error is smaller than the estimation of
  the discretization error. Using this criterion and
  an adaptive mesh refinement technique, we obtain  an efficient and robust method for
  the numerical solution of PDEs, which is demonstrated by several numerical experiments.
\keywords{Goal-oriented error estimates \and algebraic errors \and BiCG method \and adaptivity}
 \subclass{65N15 \and 65N30 \and 65F10 \and 15A06}
\end{abstract}

\section{Introduction}

Let $\Om \subset \IR^d$ be a polygonal domain with the boundary  $\gom=\partial\Om$.
We consider an abstract partial differential equation in the form
\begin{align}
  \label{eq:PP}
  \krl \u = f,
\end{align}
where $\u:\Om \to \IR$ is the unknown solution,
$\krl$ is a linear differential operator and $f$ is the right-hand side.
The problem \eqref{eq:PP} has to be accompanied by suitable boundary conditions.
Further, let $\Vh$ be a finite-dimensional space of functions, ($\dim\Vh=\Nh< \infty$),
where an approximation of $\u$ is sought. 
We note that the space $\Vh$ is not given a priori but it is generated as
a sequence of spaces by a suitable mesh adaptive method.
By $\uh\in\Vh$ we denote the discrete solution which approximates
$u$.

In many practical application, 
we are not interested
in the solution of \eqref{eq:PP} itself but rather in a {\em quantity of interest},
which is the value of a certain, a priori known,
solution-dependent {\em target functional} $\J(\u)$.
For example, the functional $\J$ is a mean value of the solution
over a subset of the computational domain or its boundary.
Therefore, we need not to estimate the  error
$\norm{\u-\uh}{}$
in an energy norm (cf. \cite{AinsworthOden,verfurth-book2}) but
the error $\J(u)-\J(\uh)$.
We assume that $J$ is a linear functional.

The necessity to estimate the error of the target functional gives rise
to the  {\em goal-oriented error estimates},
cf. the pioneering works summarized in \cite{RannacherBook,BeckerRannacher01,GileSuli02}.
This approach was further developed for many problems, let us mention
\cite{OdenPrudhomme_CAMWA01,Korotov_JCAM06,SolDemko04} dealing with
linear elliptic problems,
\cite{KuzminMoller_JCAM10} dealing with steady linear hyperbolic problems.
Moreover, extensions to nonlinear problem were presented, e.g., in 
\cite{Rey2Gosselet_CMAME14,Rey2Gosselet_CMAME15} for  elasticity
problems
and, e.g., in 
\cite{LoseilleDervieuxAlauzet_JCP10,Richter_IJNMF10,HH06:SIPG2,GeHaHo09,Hartman08,BalanWoopenMay16}
for computational fluid dynamics, 
for a survey see \cite{FidkowskiDarmofal_AIAA11}.

The goal-oriented error estimates 
require, except the solution of the original ({\em primal}) problem
\eqref{eq:PP},
also to solve the {\em dual} (or adjoint) problem
\begin{align}
  \label{eq:DP}
  \krl^* \z = \J,
\end{align}
where $\krl^*$ is the {\em dual operator} to $\krl$, $\J$ is the target functional
and $\z:\Om\to\IR$ is the dual solution.
The problem \eqref{eq:DP} has to be accompanied by suitable boundary conditions as well.

It is necessary to approximate the solution of \eqref{eq:DP} by $\zhp\in\Vhp$
where $\Vhp$ is a richer space than $\Vh$.
One possibility is to discretize \eqref{eq:DP} directly on $\Vhp$
(e.g., \cite{Korotov_JCAM06})
but then we have to solve two different algebraic problems. Another way
is to discretize \eqref{eq:DP} on $\Vh$ and define $\zhp:= \krr(\zh)$,
where $\zh\in\Vh$ is a numerical approximation of \eqref{eq:DP} and
$\krr:\Vh\to\Vhp$ is a suitable higher order reconstruction
(e.g., \cite{Richter_IJNMF10,CarpioPrietoBermejo_SISC13}).

The advantage of the latter  approach is that 
the discretizations of the primal and dual problems \eqref{eq:PP} and \eqref{eq:DP}
using the same space $\Vh$
are equivalent to two mutually transposed linear algebraic systems
which can be beneficial in practical solutions. Namely, we obtain
\begin{align}
  \label{i0}
  \mA\bx = \bb\qquad\mbox{and}\qquad \mA^{\T}\by=\bc,
\end{align}
where $\mA\in\IR^{\Nh\times\Nh}$ is the matrix arising from the discretization
of $\krl$ from~\eqref{eq:PP},  $\mA^{\T}$ is the transpose matrix of $\mA$,
$\bb\in\IR^{\Nh}$ and $\bc\in\IR^{\Nh}$ represent the discretization of the
right-hand sides of \eqref{eq:PP} and \eqref{eq:DP}, respectively, and
$\bx\in\IR^{\Nh}$ and $\by\in\IR^{\Nh}$ are the vectors corresponding to
the discrete solutions $\uh\in\Vh$ and $\zh\in\Vh$ of  \eqref{eq:PP} and \eqref{eq:DP},
respectively.

In many situations, it is advantageous to solve systems \eqref{i0} iteratively
since (i) approximate solutions satisfying the prescribed
tolerance are sufficient as an output of the computation; (ii)  very good initial
approximations of $\bx$ and $\by$ are typically available from previous level of
mesh adaptation. 
Let  $\uhk\in\Vh$ and $\zhk\in\Vh$ denote functions corresponding to the
$k$-th iterations $\bxk$ and $\byk$, respectively.
Then the computed approximations $\uhk$ and $\zhk$ are influenced also
by the algebraic error arising from the inexact solution of systems \eqref{i0}.
Let us note that even if a direct solver for \eqref{i0} is used then these systems
are not solved exactly  (only on the level of machine accuracy) and the computed
approximations suffers from algebraic errors too, see \cite{ArioliLiiesenMiedlerStrakos13}.

Having only approximate solution of \eqref{i0}, the Galerkin orthogonality of the error
is violated and many standard a posteriori error estimation techniques can not be used.
Therefore, the algebraic error has to be taken into account in a posteriori error analysis.
Additionally, the algebraic error estimation is important for the setting of
a suitable stopping criterion for iterative solvers and therefore for the
optimization of the computational costs.
The algebraic error estimates in the framework of energy norms were treated, e.g., in
\cite{Arioli_NM04,Picasso_CNME09,VohralikStrakos10} for conforming finite element method in
combination with the conjugate gradient method.
In the framework of goal-oriented error estimates,
a possible violation of Galerkin orthogonality was taken into account, e.g., 
in \cite{KuzminKorotov_MCS10}. 

The influence of the algebraic error
in goal-oriented error estimates was considered 
in \cite{Meidner2009Goal} for iterative algebraic solvers
(including multigrid methods). The computational error was expressed and
estimated as a sum of
the {discretization} and {algebraic} errors.
The primal and dual algebraic problems \eqref{i0}
were solved alternatively and after some number of iterations
(or one multigrid cycle), the algebraic and discretization components
of the errors were estimated. Then the iterative process was stopped when
the algebraic error estimate  was sufficiently smaller than the discretization error estimate.
We also mention the recent paper \cite{MallikVohralikYousef_JCAM20} using
  a completely deferent approach  for elliptic problems.
  This technique is based on $H(\mathrm{div})$-conforming flux reconstructions
  and $H^1$-conforming potential reconstructions and yields to
  a guaranteed upper error bound.

In \cite{DolRoskovec_AM17}, we employed
an approach similar to \cite{Meidner2009Goal}
and studied the convergence of the
error estimator with respect to algebraic iterations of the primal and dual problems.
We observed a delay in the convergence for insufficiently accurate resolution
of the primal or the dual problem.
We arrived to the conclusion that it is difficult to control the efficiency of
the computational process when the primal and dual problems \eqref{i0} are solved
alternatively.

Therefore, in this paper, we develop a technique 
to solve the primal and
dual problems \eqref{i0} simultaneously using the {\em bi-conjugate gradient} (BiCG) me\-thod.
At each BiCG iteration, the approximations of $\uh$ and $\zh$ are available and
the estimate of the algebraic error can be computed and the accuracy and efficiency
can be controlled. 
Motivated by \cite{StTi2011}, we proposed a stopping criterion for
the BiCG solver which is very cheap for evaluation and in contrary to
\cite{Meidner2009Goal,DolRoskovec_AM17}, it does not require the evaluation of
the discretization error estimator. The use of the BiCG solver for
the primal and dual problems \eqref{i0} with the proposed criterion
in a combination the mesh adaptation leads to an efficient
numerical method for the solution of \eqref{eq:PP}.

Since the presented approach can be applied to \eqref{eq:PP}
with a general linear operator $\krl$
discretized by any numerical scheme based on a variational formulation,
we express the problem considered and its numerical approximation
only in the abstract form.
However, in order to demonstrate the applicability of this technique,
we present the numerical solution of purely elliptic and convection-diffusion problems by
the $hp$-adaptive discontinuous Galerkin method (DGM) on possibly anisotropic meshes.

The content of the rest of the paper is the following.
In Section~\ref{sec:GO}, we shortly summarize the 
goal-oriented error estimate technique including two variants of the expression
of the algebraic errors.
In Section~\ref{sec:algeb}, we introduce the algebraic representation of the
discrete problems and discuss several ways for the evaluation of the quantity of interest
and their corresponding error estimates. These possibilities are numerically tested
in Section~\ref{sec:numerF}, where we solve the Laplace problem on fixed meshes.
Moreover, in Section~\ref{sec:adapt}, we describe the mesh adaptation process
and proposed new stopping criteria for the BiCG solver.
Their computational performance is demonstrated in Section~\ref{sec:numerA}
for the Laplace and convection-dominated problems.
We conclude with several remarks and discuss a  possible extension of this technique
for nonlinear problems in Section~\ref{sec:sum}.

\section{Framework for the goal-oriented error estimates}
\label{sec:GO}

We briefly recall a general framework for the goal-oriented error estimates.
More details can be found, e.g., in \cite{BeckerRannacher01,GileSuli02}.

\subsection{Primal problem}

Let the {\em weak formulation} of the primal problem \eqref{eq:PP} be given by
\begin{align}
  \label{eq:PPw}
  a( \u, \vp) = \ell(\vp)\qquad \forall \vp\in \V,
\end{align}
where $\u\in \V$ is a weak solution, $a(\cdot,\cdot):\V\times \V\to\IR$ is a bilinear form,
$\ell(\cdot):\V\to\IR$ is a linear form and $V$ is a Hilbert space.
We assume that \eqref{eq:PPw} is well-posed, i.e., it admits a unique weak solution.

For the numerical approximation of \eqref{eq:PPw}, let $\Vh,\ h>0$ be
a finite element space of functions defined on $\Om$, typically piecewise
polynomial functions related  to the partition of $\Om$ onto a set of finite elements $\Th$.
Moreover, let $\Wh$ be a functional space such that $\V\subset \Wh$ and
$\Vh\subset\Wh$. For conforming finite element methods ($\Vh\subset\V$), we simply
put $\Wh:=\V$. However, for  nonconforming methods (the case $\Vh\not\subset\V$),
the choice of $\Wh$ is more
delicate. E.g., for {\em discontinuous Galerkin method}, we employ the so-called
{\em broken Sobolev space} and put
\begin{align}
  \label{HsTh}
  \Wh:= H^s(\Th) = \{ v \in L^{2}(\Om); \, \restr{v}{K} \in H^{s}(K) \,
  \forall K \in \Th\},
\end{align}
where $\Th$ is a mesh partition of $\Om$ and $s>0$ is a suitable Sobolev index,
e.g., for the second-order operator  $\krl$, we have $s=2$.

Let
\begin{align}
  \label{ah}
  \ah:\Wh\times\Wh\to\IR\qquad \mbox{and} \qquad  \ell_h:\Wh\to\IR
\end{align}
be a bilinear and a linear forms
corresponding to the discretization of the left-hand and right-hand sides of
\eqref{eq:PPw}, respectively,  by a particular numerical method.

We say that $\uh\in\Vh$ is the {\em discrete solution} of
the {\em primal problem} \eqref{eq:PPw} if
\begin{align}
  \label{eq:PPh}
    \ah(\uh,\vp_h) = \ell_h(\vp_h)\qquad \forall \vp_h\in\Vh.
\end{align}
We assume that the numerical scheme \eqref{eq:PPh} is {\em consistent}, i.e.,
\begin{align}
  \label{eq:cons}
  \ah(\u,\vp) = \ell_h(\vp)\qquad \forall \vp\in\Wh
\end{align}
where $\u\in\V$ is the weak solution of \eqref{eq:PPw}.
This implies the {\em Galerkin orthogonality of the error} of the primal
problem
\begin{align}
  \label{eq:GO}
  \ah(\uh - u, \vp_h) = 0 \qquad\forall \vp_h\in\Vh.
\end{align}
Finally, we define the {\em residual of the primal problem} by
\begin{align}
  \label{res}
  \res{\uh}{\vp}:= \ell_h(\vp) - \ah(\uh, \vp) = \ah(\u-\uh,\vp) ,\quad \vp\in\Wh,
\end{align}
where the last equality follows from the consistency \eqref{eq:cons}.


\subsection{Quantity of interest and the dual problem}

As mentioned in the introduction, we are interested in
 a {\em sufficiently accurate} approximation of
the {\em quantity of interest} $J(u)\in\IR$, where
\begin{align}
  \label{Ju}
  \J:\Wh\to \IR
\end{align}
is  a linear functional. 
Typically, it is defined as a weighted mean value of $u$ over the computational
domain $\Om$ or its boundary $\gom$.

In order to estimate the error $\J(\u) - J(\uh)$, we consider the 
{\em adjoint} (or {\em dual}) {\em problem} \eqref{eq:DP} 
and its discretization. 
We say that $\zh$ is the {\em discrete solution} of the {\em dual problem}
\eqref{eq:DP} if
\begin{align}
  \label{eq:DPh}
  \ah(\psi_h, \zh) = \J(\psi_h) \qquad \forall \psi_h\in\Vh,
\end{align}
where $\ah$ and $\J$ are given by \eqref{ah} and \eqref{Ju},
respectively.
Moreover, we assume that the numerical scheme \eqref{eq:DPh} is
{\em adjoint consistent}, i.e.,
\begin{align}
  \label{eq:consD}
  \ah(\psi, \z) = \J(\psi)\qquad \forall \psi\in\Wh,
\end{align}
where $\z$ is the weak solution of the dual problem \eqref{eq:DP}.
This implies the {\em Galerkin orthogonality of the error} of the dual problem
\begin{align}
  \label{eq:GOD}
  \ah(\psi_h, \zh - \z) = 0 \qquad\forall \psi_h\in\Vh.
\end{align}
Finally, we define the {\em residual of the dual problem} by
\begin{align}
  \label{resD}
  \resD{\zh}{\psi}:= \J(\psi) - \ah(\psi,\zh) = \ah(\psi, \z-\zh),\quad \psi\in\Wh,
\end{align}
where the last equality follows from the adjoint consistency \eqref{eq:consD}.


The problems \eqref{eq:PPh} and \eqref{eq:DPh} represent two linear algebraic systems
whose efficient solution by an iterative solver is developed in \S\ref{sec:algeb}.
Due to iterative and rounding errors, the ``exact'' discrete solutions $\uh$ 
and $\zh$ are not available, we have only
their approximation $\uhk\in\Vh$ and $\zhk\in\Vh$, $k=1,2,\dots$.
However, they do not fulfil
the Galerkin orthogonalities \eqref{eq:GO} and \eqref{eq:GOD}.

\subsection{Abstract goal-oriented error estimates}
\label{sec:GOEE}

Obviously, for adjoint consistent discretization we have (due to
\eqref{eq:cons} and \eqref{eq:consD}) the equivalence
between the quantity of interest $\J(\u)$ and the right-hand side of \eqref{eq:PPw}
evaluated for the dual solution $\z,$ i.e. 
\begin{align}
  \label{eq:eqi}
 \ell_h(\z) = \ah(\u, \z) = \J(\u).
\end{align}
Similarly, from \eqref{eq:PPh} and  \eqref{eq:DPh}, we obtain the discrete variant
of \eqref{eq:eqi} as
\begin{align}
  \label{eq:eqi_h}
 \ell_h(\zh) = \ah(\uh, \zh) = \J(\uh),
\end{align}
where $\uh$ and $\zh$ are the discrete solutions of the primal and dual problems,
respectively.
Consequently, we have an {\em error equivalence}
\begin{align}
  \label{eq:eqi_err}
  \J(\u -\uh) =  \J(\u) - \J(\uh) = \ell_h(\z) - \ell_h(\zh) = \ell_h(\z -\zh).
\end{align}
Therefore, the difference $\ell_h(\z -\zh)$ can be used as an error estimate
of the quantity of interest as well.

First, we present the primal and dual error identities for the error of the
quantity of interest for the algebraically exact  discrete solution $\uh$
and $\zh$ of \eqref{eq:PPh} and \eqref{eq:DPh}, respectively.
Using  the adjoint consistency \eqref{eq:consD}, 
the Galerkin orthogonality \eqref{eq:GO}
and 
relation  \eqref{res},
we get the {\em primal error identity} 
\begin{align} 
 \label{eq:JP}
 \J(\u - \uh) &=  \ah(\u - \uh, \z)=  \ah(\u - \uh, \z - v_h)  \\
   &= \res{\uh}{\z -  v_h}  \quad \forall v_h \in \Vh, \notag
\end{align}
where $\res{\zh}{\cdot}$ denotes the residual of the primal problem given by \eqref{res}.

Similarly, exploiting in addition the Galerkin orthogonality \eqref{eq:GOD} and
relation  \eqref{resD},
we get the {\em dual error identity}
\begin{align} 
\label{eq:JD}
\J(\u - \uh) &=  \ah(\u - \uh, \z - \zh)
=   \ah(\u - w_h, \z - \zh ) \\
&=   \resD{\zh}{\u -  w_h}  \quad \forall w_h \in \Vh, \notag
\end{align}
where $\resD{\zh}{\cdot}$ denotes the residual of the dual problem given by \eqref{resD}.

\subsection{Abstract goal-oriented error estimates including algebraic errors}
\label{sec:GOEE_alg}

As mentioned above, the discrete solutions $\uh$ and $\zh$ fulfilling
\eqref{eq:PPh} and \eqref{eq:DPh}, respectively, are non-available but
we have only $\uhk\in\Vh$  and $\zhk\in\Vh$, $k=1,2,\dots$  denoting their approximations
given by an algebraic iterative solver. Hence, the error representations
\eqref{eq:JP} and \eqref{eq:JD} are useless and we have to estimate
the error $J(\u-\uhk)$.

In the same way as in \cite{Meidner2009Goal,DolRoskovec_AM17},
we employ the adjoint consistency \eqref{eq:consD}, identity $\z = \z-v_h + v_h,\ v_h\in\Vh$
and 
relation  \eqref{res},
we get the {\em primal error identity} including {\em algebraic errors}
\begin{align} 
 \label{eq:JP_alg0}
 \J(\u - \uhk) &=  \ah(\u - \uhk, \z)=  \ah(\u - \uhk, \z - v_h) + \ah(\u - \uhk, v_h)
 \\
 &   = \res{\uhk}{\z -  v_h} + \res{\uhk}{v_h} \quad \forall v_h \in \Vh,\notag
\end{align}
where $\res{\uhk}{\cdot}$ 
is given by \eqref{res}.
Inserting $v_h:=\zhk$ in \eqref{eq:JP_alg0}, we obtain
\begin{align} 
 \label{eq:JP_alg}
 \J(\u - \uhk) = \res{\uhk}{\z -  \zhk} + \res{\uhk}{\zhk}=: \errSk + \errAk,
\end{align}
where the quantity $\errSk=\res{\uhk}{\z -  \zhk}$
represents the discretization
error of the primal problem since it coincides with \eqref{eq:JP} for $\uhk=\uh$.
Further, the quantity $\errAk=\res{\uhk}{\zhk}$
represents the algebraic
error of the primal problem since it vanishes for  $\uhk=\uh$.
Let us note that \eqref{eq:JP_alg0} and \eqref{eq:JP_alg} hold
without the Galerkin orthogonality \eqref{eq:GO}.

In order to derive the analogue of \eqref{eq:JD} including algebraic error,
we take into account the error equivalence \eqref{eq:eqi_err}.
The quantity $\ell_h(\z - \zhk)$ exhibits the analogue of
the error $ \J(\u - \uhk)$ since both quantities are equal for $\uhk=\uh$ and $\zhk=\zh$.
However, if $\uhk$ and $\zhk$ are arbitrary iterations then,
in general,
\begin{align}
  \label{eq:equi_alg}
  \J(\u - \uhk) \not= \ell_h(\z - \zhk).
\end{align}
Nevertheless,  we show in Section~\ref{sec:algeb} (cf. \eqref{xi8}) that
if $\uhk$ and $\zhk$ are obtained by the bi-conjugate gradient (BiCG) method
with vanishing initial approximations 
(in the exact arithmetic) then 
\begin{align}
  \label{eq:equi_algB}
  \J(\u - \uhk) =\ell_h(\z - \zhk)\qquad \mbox{ for } k=1,2,\dots.
\end{align}

Now, using the consistency \eqref{eq:consD}, identity $\z = \z-w_h + w_h,\ w_h\in\Vh$
and 
relation  \eqref{res},
we get the {\em dual error identity} including {\em algebraic errors}
\begin{align} 
  \label{eq:JD_alg0}
  \ell_h(\z - \zhk) &=  \ah(\u , \z-\zhk)=  \ah(\u - w_h, \z - \zhk) + \ah( w_h, \z - \zhk)
  \\
  &   = \resD{\zhk}{\u -  w_h} + \resD{\zhk}{w_h} 
  \quad \forall w_h \in \Vh,\notag
\end{align}
where $\resD{\zhk}{\cdot}$ 
is given by \eqref{resD}.
Putting $w_h:=\uhk$ in \eqref{eq:JD_alg0}, we obtain
\begin{align} 
  \label{eq:JD_alg}
  \ell_h(\z - \zhk) = \resD{\zhk}{\u -  \uhk} + \resD{\zhk}{\uhk} =: \errDSk + \errDAk,
\end{align}
where, similarly as the quantities $\errSk$ and $\errAk$ in \eqref{eq:JP_alg},
the quantities $\errDSk$ and $\errDAk$ represent the discretization and algebraic
error of the dual problem.
Let us note that \eqref{eq:JD_alg0} and  \eqref{eq:JD_alg} hold
without the Galerkin orthogonality \eqref{eq:GOD}.

\subsection{Computable goal-oriented error estimates}
\label{sec:GOEE_com}

Whereas the algebraic errors $\errAk$ and $\errDAk$ from
\eqref{eq:JP_alg} and \eqref{eq:JD_alg}, respectively,
are computable quantities,
the discretization errors $\errSk$ and $\errDSk$ 
require the knowledge of the exact dual and primal and solutions $\z$ and $\u$,
respectively.
In practical computations, they should be approximated by a higher-order
reconstruction denoted here by
\begin{align}
  \label{krr}
  \huh=\krr(\uhk),\qquad \hzh=\krr(\zhk),
\end{align}
where  $\uhk\in\Vh$  and $\zhk\in\Vh$ denote approximations of $\uh$ and $\zh$,
$\krr:\Vh\to\Vhp$  denotes a reconstruction
operator and $\Vhp$ is a ``richer'' finite dimensional space,
for examples, see the papers cited in Introduction.

Replacing $\u$ and $\z$ in \eqref{eq:JP_alg} and \eqref{eq:JD_alg} by
$\huh$ and $\hzh$, we obtain the computable approximations $\etaSk$ and $\etaDSk$
of the discretization errors $\errSk$ and $\errDSk$ by
\begin{align}
  \label{etas}
  &\errSk = \res{\uhk}{\z -  \zhk} \approx \res{\uhk}{\hzh -  \zhk} =: \etaSk, \\
  &\errDSk = \resD{\zhk}{\u -  \uhk} \approx \res{\zhk}{\huh -  \uhk} =: \etaDSk,\qquad
  k=1,2,\dots.
  \notag
\end{align}
The quantities $\etaSk$ and $\etaDSk$ are called the estimates of the
primal and dual discretization errors, respectively.

The approximations \eqref{etas} together with \eqref{eq:JP_alg} and \eqref{eq:JD_alg}
lead to the computable estimate of the total error by
\begin{align}
  \label{EEP}
  \J(\u - \uhk) & \approx   \etaSk + \etaAk, \qquad 
  \ell_h(\z - \zhk)  \approx \etaDSk + \etaDAk, \quad k=1,2,\dots,
\end{align}
where, in order to have a consistent notation, we put
\begin{align}
  \label{etasA}
  \etaAk := \errAk = \res{\uhk}{\zhk}, \quad
  \etaDAk := \errDAk = \resD{\zhk}{\uhk}.
\end{align}
The terms $\etaAk$ and $\etaDAk$ are equal to the
{\em algebraic errors} of
the primal and dual problems \eqref{eq:PPh} and \eqref{eq:DPh}, respectively,
and they are independent on the used
higher-order reconstruction $\krr$ from \eqref{krr}.

\begin{remark}
  \label{rem:guaranteed}
  The approximations \eqref{EEP} do not give a guaranteed error estimate
  since, as usual,  we neglected the terms
  $\res{\uhk}{\z -  \hzh} $ and $\res{\zhk}{\u -  \huh} $.
  Hence, the guaranteed error estimate requires an estimation of these terms.
  We refer to \cite{NochettoVeeserVerani_IMAJNA09,Ainsworth2012Guaranteed}
  where this problem was treated for symmetric elliptic problems discretized by
  the conforming finite element method. This is also a subject of our further research.
\end{remark}

\subsection{An alternative representation of the algebraic error}
\label{sec:alter}

In \eqref{eq:JP_alg}, we expressed the discretization and algebraic parts
of the computational error using the approach from \cite{Meidner2009Goal},
namely 
\begin{align} 
 \label{JP1}
 \J(\u - \uhk) = \res{\uhk}{\z -  \zhk} + \res{\uhk}{\zhk} = \errSk +\errAk.
\end{align}
On the other hand,
it is possible to decompose the total error $\J(\u - \uhk)$
in the alternative way as
\begin{align}
  \label{JP2}
  \J(\u - \uhk) = \J(\u - \uh) + \J(\uh - \uhk)=: \estSk + \estAk.
\end{align}
The first term in \eqref{JP2} is the discretization error given by
\eqref{eq:JP} and it satisfies
\begin{align}
  \label{JP3}
   \estSk = \J(\u - \uh) =  \res{\uh}{\z -  v_h}, \quad  v_h \in \Vh
\end{align}
and it is independent of $k$. Obviously, although $\errSk$ in  \eqref{JP1}
and $\estSk$ in \eqref{JP3} 
both represent the discretization error, they  differ.
Moreover, $\estSk$ cannot be evaluated even if we replace $\z$ by $\hzh$ since
$\uh$ is unavailable.

The second term $\estAk$ in \eqref{JP2}  represents the algebraic error.
It cannot be expressed in a residual form but
in Section~\ref{sec:algeb}, we present a technique
which is able to estimate this term in a cheap way.

Hence, we have two representations of the algebraic error, the first one from \eqref{JP1}
given by $\etaAk=\res{\uhk}{\zhk}$ and the second one from \eqref{JP2} given by
$\estAk:=\J(\uh - \uhk)$.
Both representations are ``algebraically consistent'' which means that
if $\uhk\to \uh$ and $\zhk\to \zh$ for $k\to \infty$  then
\begin{align}
  \label{JP4}
  \etaAk \to 0 \quad \mbox{ for } k\to \infty
  \qquad \mbox{ and }\qquad 
  \estAk \to 0 \quad \mbox{ for } k\to \infty.
\end{align}
However, the speed of convergence for $\etaAk$ and $\estAk$ can differ substantially,
  see numerical experiments in Section~\ref{sec:numerA}.

Moreover, for $\uhk\to \uh$ and $\zhk\to \zh$, both discretization errors
$\errSk$ and $\estSk$ are closer and closer
since (cf. \eqref{JP1} and \eqref{JP2})
\begin{align}
  \label{JP4a}
  | \errSk - \estSk| = | \estAk - \errAk| \leq | \estAk| + |\errAk| \to 0.
\end{align}

Finally, we present the dual analogue of \eqref{JP2}--\eqref{JP3} in the form
\begin{align}
  \label{JP5}
  \ell_h(\z - \zhk) &= \ell_h(\z - \zh) + \ell_h(\zh - \zhk) 
  =: \estDSk + \estDAk,
\end{align}
where the term $\estDAk$ can be estimated in the same manner as $\estAk$ in \eqref{JP2},
cf. Section~\ref{sec:algeb}.

\section{Solution of primal and dual discretized problems}
\label{sec:algeb}

In this section, we introduce the algebraic representation of
the goal-oriented error estimates from the previous section, present the BiCG method
allowing a simultaneous solution of the primal and dual problems and discuss
several possibilities  algebraic errors estimates and
stopping criteria for the iterative solver.

\subsection{Algebraic representation}
\label{sec:alg1}

Let $\{\vp_i:\Om\to\IR,\ i=1,\dots, \Nh \}$ be a basis of the finite-dimensional space $\Vh$.
We define the matrix $\mA\in\IR^{\Nh\times\Nh}$ by
\begin{align}
  \label{mA}
  \mA = \{\mA_{i,j}\}_{i,j=1}^{\Nh},\quad \mA_{i,j} := \ah(\vp_i,\vp_j),\ i,j=1,\dots, \Nh,
\end{align}
where $\ah$ is the bilinear form defined by \eqref{ah}.
Then the primal and dual discrete problems \eqref{eq:PPh} and \eqref{eq:DPh}
are equivalent to the solution of
two linear algebraic systems
\begin{align}
  \label{AA1}
  \mA\bx = \bb\qquad\mbox{and}\qquad \mA^{\T}\by=\bc,
\end{align}
respectively, where 
$\bx\in\IR^\Nh$ and  $\by\in\IR^\Nh$ are the algebraic representation
of the primal and dual solutions given by
\begin{align}
  \label{AA2}
  & \bx = (\x_1,\dots, \x_{\Nh})^{\T} \ \leftrightarrow \uh =\sum\nolimits_{i=1}^\Nh \x_i \vp_i, \\
  \quad 
  \mbox{and} \quad &   \by = (\y_1,\dots, \y_{\Nh})^{\T} \ \leftrightarrow
  \zh =\sum\nolimits_{i=1}^\Nh \y_i \vp_i,
  \notag
\end{align}
respectively, and $\bb\in\IR^\Nh$ and  $\bc\in\IR^\Nh$ are the algebraic representation
of the right-hand sides of the primal and dual problems given by
\begin{align}
  \label{AA3}
  & \bb = (\b_1,\dots, \b_{\Nh})^{\T} \ \leftrightarrow \b_i = \ell_h(\vp_i),\ i=1,\dots,\Nh,
  \\ 
   \mbox{and} \quad
   &  \bc = (\c_1,\dots, \c_{\Nh})^{\T} \ \leftrightarrow \c_i =\J(\vp_i),\quad i=1,\dots,\Nh,
   \notag
\end{align}
respectively.
Using \eqref{mA}--\eqref{AA3}, we obtain the equivalencies
\begin{align}
  \label{AA4}
  \J(\uh)  & 
  = \sum_{i=1}^\Nh \x_i \J(\vp_i)
  = \bcT \bx = \bcT \mA^{-1} \bb 
   = \by^{\T} \bb 
  = \sum_{i=1}^{\Nh} \y_i \ell_h(\vp_i)
  = \ell_h(\zh), 
\end{align}
which exhibits the algebraic analogue of \eqref{eq:eqi_h}.

Similarly, as in \eqref{AA2} let $\bxk$ and $\byk$ be the algebraic analogues
of the approximations $\uhk$ and $\zhk$, respectively.
The corresponding residuals of \eqref{AA1} are given by 
\begin{align}
  \label{AA5}
  \brk := \bb - \mA \bxk =\mA (\bx - \bxk) \quad \mbox{ and } \quad
  \bsk := \bc - \mA^{\T} \byk = \mA^{\T} (\by -  \byk).
\end{align}
By a standard manipulation,
we derive the following correspondence between the discretization and its algebraic
representation which will be used in next paragraphs. Similarly as in \eqref{AA4},
we have
\begin{align}
  \label{aa1}
  \J(\uhk)  & = \bcT \bxk, \qquad 
  \ell_h(\zhk)   = \bykT \bb.
\end{align}
From \eqref{mA} and \eqref{AA2}, we obtain
\begin{align}
  \label{aa3}
  \ah(\uhk, \zhk ) & = \bykT \mA \bxk.
\end{align}
Additionally, employing \eqref{AA1} and \eqref{AA5}, we derive
\begin{align}
  \label{aa4}
  \ah(\uh - \uhk, \zh - \zhk ) & = (\by - \byk)^{\T} \mA (\bx - \bxk) = 
  \bskT \mA^{-1} \brk.
\end{align}
Finally, using \eqref{res}, \eqref{resD}, \eqref{etasA} and \eqref{AA5}--\eqref{aa3},
we obtain
\begin{align}
  \label{aa5}
  \etaAk & = \ell(\zhk) - \ah(\uhk ,\zhk)
  = \bykT \bb - \bykT \mA \bxk
  = \bykT \brk, \\
  \etaDAk & = \J(\uhk) - \ah(\uhk ,\zhk)
  = \bcT \bxk - \bykT \mA \bxk
  = \bskT \bxk. \notag
\end{align}

\subsection{Approximating the quantity of interest using iterates $\bxk$ and $\byk$}
\label{sec:eqi}

Let the approximations $\bxk$ and $\byk$, and the corresponding residual vectors $\brk$ and $\bsk$, computed by some iterative method for solving linear systems \eqref{AA1}, be given.
We introduce several possibilities, how to approximate the quantity of interest
$\J(\uh) = \bcT \mA^{-1} \bb = \ell_h(\zh)$ using these vectors.
For each possibility, we express the quantity of interest as a sum of two terms,
where the former one is a computable value approximating
the quantity of interest and the latter one is an incomputable term which represents
the error of the approximation. Note that the identities derived below rely only on the relations
\begin{equation}\label{eqn:residuals}
  \brk = \bb - \mA \bxk \quad \mbox{ and } \quad
  \bsk = \bc - \mA^{\T} \byk.
\end{equation}

\begin{itemize}
\item [\PPa] Using \eqref{AA5}, a simple manipulation gives
  \begin{align}
    \label{BB1}
    \bcT \mA^{-1} \bb = \bcT\bx = \bcT\bxk +  \bcT( \bx - \bxk)
    =  \underbrace{\bcT\bxk}_{\mathrm{approximation}}
    + \underbrace{\bcT \mA^{-1} \brk}_{\mathrm{error}}.
  \end{align}
  The first term on the right-hand side of \eqref{BB1} is computable
  and then it can be used for the approximation
  $\bcT \mA^{-1} \bb \approx \bcT\bxk$ and the second term
  $\bcT \mA^{-1} \brk$ represents the corresponding error.
  Similarly, for the dual form, we have
  \begin{align}
    \label{BB2}
    \bcT \mA^{-1} \bb = \by^{\T}\bb = \bykT \bb +  ( \by - \byk)^{\T}\bb
    = \underbrace{\bykT \bb}_{\mathrm{approximation}} +
    \underbrace{ \bskT\mA^{-1} \bb}_{\mathrm{error}},
  \end{align}
  thus $\bykT \bb$ is a computable approximation of $\bcT \mA^{-1} \bb$
  and $\bskT\mA^{-1} \bb$ the corresponding error.
\item[\PPb] We carry out more sophisticated manipulation, take into account
  \eqref{AA1}, \eqref{AA5}, \eqref{aa4} and  get
   \begin{align}
    \label{BB3}
    \bcT \mA^{-1} \bb & = \bcT\bxk +  \bcT \mA^{-1}\mA( \bx - \bxk) 
     =  \bcT\bxk + \by^{\T} \mA(\bx-\bxk)   \\
    & =   \bcT\bxk + \bykT (\bb- \mA\bxk) + (\by-\byk)^{\T}\mA (\bx- \bxk)
    \notag \\
    & = \underbrace{\bcT\bxk + \bykT \brk}_{\mathrm{approximation}}
    + \underbrace{\bskT\mA^{-1} \brk}_{\mathrm{error}}. \notag
   \end{align}
   Using the algebraic-discretization equivalence relations \eqref{AA4},
   \eqref{aa1}, \eqref{aa4} and \eqref{aa5}, the identity \eqref{BB3} can be written in the equivalent
   form
   \begin{align}
     \label{BB3a}
     \J(\uh) = \underbrace{\J(\uhk)+ \etaAk}_{\mathrm{approximation}}
     + \underbrace{\ah(\uh - \uhk, \zh - \zhk )}_{\mathrm{error}}.
   \end{align}
   Similarly, we can derive the dual relation
   \begin{align}
     \label{BB4}
     \bcT \mA^{-1} \bb
     = \underbrace{\bykT\bb + \bskT \bxk}_{\mathrm{approximation}}
     + \underbrace{\bskT\mA^{-1} \brk}_{\mathrm{error}},
   \end{align}
   which is equivalent to
   \begin{align}
     \label{BB4a}
     \ell_h(\zh) =\underbrace{ \ell_h(\zhk) + \etaDAk}_{\mathrm{approximation}}
     + \underbrace{\ah(\uh - \uhk, \zh - \zhk )}_{\mathrm{error}}.
   \end{align}
\end{itemize}
As mentioned at the beginning of this section, both evaluations {\PPa} and {\PPb} can be used for any iterative solver
generating approximations $\bxk$ and $\byk$, and residuals $\brk$ and $\bsk$, $k=1,2,\dots$.
If the norms of the residual vectors $\brk$ and $\bsk$ tend to be small, 
one can expect that
\begin{align}
  \label{BB5}
  | \bskT\mA^{-1} \brk | \ll |\bcT \mA^{-1} \brk| \qquad \mbox{and}\qquad
  |\bskT\mA^{-1} \brk| \ll |  \bskT\mA^{-1} \bb|.
\end{align}
In other words, one can expect that the approximation  {\PPb} is more accurate than
{\PPa} since the corresponding error terms tend to be smaller.

\subsection{Approximating the quantity of interest using BiCG iterates}
\label{sec:bicg}

As mentioned in Introduction, our aim is to solve systems \eqref{AA1} by an iterative method, which allows to solve the primal
and dual problems simultaneously. Since we intend to solve large and sparse systems, we need to pick a method with low memory requirements.
Then, a natural choice is to use
the {\em preconditioned bi-conjugate gradient} (BiCG) method (Algorithm~\ref{alg:BiCG}) introduced in \cite{Fl1976}; see also \cite{B:BaBeCh1994,StTi2011}.
Let $\mP$ be a suitable preconditioner (its choice depends on the particular discretization scheme).
 BiCG is a short-term recurrence Krylov subspace method which generates (if no breakdown occurs) approximations $\bxk \in \bx_0 + \mP^{-1}\mathcal{K}_k (\mA \mP^{-1},\br_0)$ and $\byk
\in \by_0 + \mP^{-\T}\mathcal{K}_k (\mA^{\T}\mP^{-\T},\bs_0)$ such that 
\begin{equation}\label{eqn:bicg}
    \brk \perp \mP^{-\T}\mathcal{K}_k (\mA^{\T}\mP^{-\T},\bs_0),\qquad \bsk \perp \mP^{-1}\mathcal{K}_k (\mA\mP^{-1},\br_0),
\end{equation}
where $\mathcal{K}_k (\mA,\br_0)$ denotes the $k$th Krylov subspace generated by $\mA$ and $\br_0$. 
The determining conditions \eqref{eqn:bicg} imply that
  \begin{align}
    \label{orthog}
    (\byk-\by_0)^{\T} \brk = 0, \qquad 
    \bskT (\bxk-\bx_0) = 0.
  \end{align}
Note that during the BiCG finite precision computations, the bi-ortho\-gonality conditions \eqref{eqn:bicg} 
are usually not satisfied, and this fact cannot be ignored in our considerations. However, not all properties of BiCG vectors are lost 
in finite precision arithmetic. Because of the choice of coefficients 
$\alpha_{k}$ and $\beta_{k+1}$ (cf. Algorithm~\ref{alg:BiCG}), 
the bi-orthogonality of two consecutive vectors (local bi-orthogonality)
is usually well preserved. This fact can be exploited to derive a more efficient way of approximating the quantities of interrest; for more details, see \cite{StTi2011}.

Algorithm~\ref{alg:BiCG} shows the BiCG algorithm which generates
sequences of primal and dual approximations $\{\bxk\}$ and  $\{\byk\}$,
respectively. 
Note that $\tilde{\br}_{k+1}=\mP^{-1} \br_{k+1}$ and $\tilde{\bs}_{k+1}=\mP^{-\T} \bs_{k+1}$
is equivalent to solving the system  
\begin{align}
  \label{mP}
  \mP\:\tilde{\br}_{k+1}=\br_{k+1}\quad\mbox{and}\quad \mP^{\T}\:\tilde{\bs}_{k+1}=\bs_{k+1},
\end{align}
respectively.
At lines 2 and 14 of Algorithm~\ref{alg:BiCG}
we compute an additional
sequence $\{\bxi_k^B\}$. The meaning of quantities $\bxi_k^B$ is explained in the text below.
The algorithm has to be furnished  by a suitable stopping criterion,
which is discussed in Section~\ref{sec:estim}.

\begin{algorithm}[ht]
  \caption{Preconditioned BiCG and approximating $\bs_{0}^{\T}\mA^{-1}\br_{0}$}
  \label{alg:BiCG}
  
\begin{algorithmic}[1]
{\large
\STATE \textbf{input} $\mA$, $\bx_{0}$, $\by_{0}$, $\mP$

\STATE $\bxi_{0}^{B}=0$

\STATE $\br_{0}=\bb-\mA\bx_{0},$ $\bs_{0}=\bc-\mA^{\T}\by_{0},$

\STATE $\bp_{0}=\mP^{-1} \br_{0}$

\STATE $\bq_{0}=\mP^{-\T} \bs_{0}$

\STATE $\tilde{\br}_{0}=\bp_{0}$ 

\FOR{$k=0,1,\dots$} 

\STATE $\alpha_{k}=\frac{\bs_{k}^{\T}\tilde{\br}_{k}}{\bq_{k}^{\T}\mA\bp_{k}}$

\STATE $\bx_{k+1}=\bx_{k}+\alpha_{k}\bp_{k}$,
\hspace{14mm}
$\by_{k+1}=\by_{k}+\alpha_k\bq_{k}$

\STATE $\br_{k+1}=\br_{k}-\alpha_{k}\mA\bp_{k}$,
\hspace{12mm}
$\bs_{k+1}=\bs_{k}-\alpha_k\mA^{\T}\bq_{k}$

\STATE $\tilde{\br}_{k+1}=\mP^{-1} \br_{k+1}$,
\hspace{18mm}
$\tilde{\bs}_{k+1}=\mP^{-\T} \bs_{k+1}$

\STATE $\beta_{k+1}=\frac{\bs_{k+1}^{\T}\tilde{\br}_{k+1}}{\bs_{k}^{\T}\tilde{\br}_{k}}$

\STATE $\bp_{k+1}=\tilde{\br}_{k+1}+\beta_{k+1}\bp_{k}$,
\hspace{8mm}
$\bq_{k+1}=\tilde{\bs}_{k+1}+{\beta}_{k+1}\bq_{k}$

\STATE $\bxi_{k+1}^{B}=\bxi_{k}^{B}+\alpha_{k}\,\bs_{k}^{\T}\tilde{\br}_{k}$

\ENDFOR
}
\end{algorithmic}
\end{algorithm}

To approximate $\bcT \mA^{-1} \bb$ in BiCG, we can use  technique from~\cite{StTi2011}, 
which is based on the local bi-orthogonality of the BiCG vectors.

First, we present some manipulations.
Let $\bx_0\in\IR^\Nh$ and $\by_0\in\IR^{\Nh}$ be the initial approximations
of $\bx$ and $\by$, respectively. Using \eqref{AA4}, \eqref{BB3}, \eqref{BB4} with $k=0$, and denoting 
\begin{align}
  \label{xis}
  \bxiP:= \bcT \bx_0 +   \by_0^{\T} \br_0, \qquad
  \bxiD:= \by_0^{\T} \bb +   \bx_0^{\T} \bs_0, \qquad 
  \bxiB:=  \bs_0 ^{\T} \mA^{-1} \br_0,
\end{align}
we obtain
the identities
\begin{align}
  \label{BA4}
    \bcT\mA^{-1}\bb = \bxiP + \bxiB = \bxiD + \bxiB\qquad \Rightarrow \quad \bxiP = \bxiD.
\end{align}

Based on technique from \cite{StTi2011} 
for estimating $\bxiB$ using $\bxiB_k$,
we are now ready to introduce the next variant of approximating
the quantity of interest.
\begin{itemize}
\item [\PPc]
  In \cite[(3.13)]{StTi2011} , it has been shown that
  \begin{align}
    \label{xi5}
    \bcT\mA^{-1}\bb  =  \underbrace{\bxiP + \bxi_{k}^{B}}_{\mathrm{approximation}}
    + \underbrace{\bskT \mA^{-1} \brk}_{\mathrm{error}},
    \quad k=0,1,\dots,
  \end{align}
with  
  \begin{align}
    \label{xi3}
     \bxi_{k}^{B}:=\sum_{n=0}^{k-1} \alpha_{n}\,\bs_{n}^{\T}\tilde{\br}_{n},
  \end{align}
  where scalars $\alpha_n$ and vectors $\bs_{n}$, $\tilde{\br}_{n}$, $n=0,1,\dots$
  are defined by Algorithm~\ref{alg:BiCG}.
  The authors of \cite{StTi2011} derive this formula using 
  the assumptions \eqref{eqn:residuals} and
  the local bi-orthogonality conditions only. Therefore, if the consecutive vectors are almost bi-orthogonal
  and if the recursively computed residuals approximately agree with the true residuals during finite precision computations,  
  then the identity \eqref{xi5} holds (up to some small inaccuracy), and   can be used for approximating
  the quantity of interest. The identity \eqref{xi5} is mathematically equivalent to the identities \eqref{BB3} and \eqref{BB4}, and represents the third possibility  of approximating the quantity of interest. The advantage of using  
  \eqref{xi5} is that we do not have to compute additional scalar products.   The quantity $\bxi_{k}^{B}$ can be computed in BiCG almost for free, since the scalar products $\bs_{n}^{\T}\tilde{\br}_{n}$ are used in BiCG to compute the coefficients $\alpha_n$ and $\beta_{n+1}$.
  
%
%
\end{itemize}

Note that in exact arithmetic, the orthogonality relations \eqref{orthog} hold. Then, comparing the identities \eqref{BB1}, \eqref{BB3}, and \eqref{xi5}, 
we obtain
$$
{\bcT\bxk + \bykT \brk}= 
{\bcT\bxk + (\byk-\by_0)^{\T} \brk+ \by_0^{\T} \brk } =  
{\bcT\bxk + \by_0^{\T} \brk }
$$
so that
\begin{align}
  \label{xi6}
  \underbrace{\bxiP + \bxi_{k}^{B}}_{\PPc}
  = \underbrace{\bcT\bxk + \bykT \brk}_{\PPb} = 
  \underbrace{\bcT\bxk}_{\PPa} + \by_0^{\T} \brk = \J(\uhk) + \by_0^{\T} \brk.
\end{align}
In particular, if $\by_0=0$, then all the evaluations {\PPa}--{\PPc} are identical for the BiCG
method in the exact arithmetic.
Similarly, from \eqref{BB2}, \eqref{BB4}, \eqref{xi5} and
using the orthogonality \eqref{orthog}, we obtain the dual counterpart relation
\begin{align}
  \label{xi7}
  \underbrace{\bxiD + \bxi_{k}^{B}}_{\PPc}
  = \underbrace{\bykT\bb + \bskT \bxk}_{\PPb}
  = \underbrace{\bykT\bb}_{\PPa} + \bskT \bx_0 = \ell_h(\zhk) + \bskT \bx_0.
\end{align}

In finite precision arithmetic, the first equalities in \eqref{xi6} and \eqref{xi7} still hold up to some small inaccuracy.
However, if the orthogonality conditions \eqref{orthog} are not (approximately) satisfied
during finite precision computations, then the second 
equalities in \eqref{xi6} and \eqref{xi7} do not (approximately) hold; for more details and examples see \cite{StTi2011}.
Note that in our experiments in Section~\ref{sec:numerF}, the orthogonality 
conditions \eqref{orthog} are well preserved and, therefore, 
the evaluations {\PPa}--{\PPc} provide almost the same results.

Finally, let us mention that if $\bx_0=0$ and $\by_0=0$, the relations  \eqref{xi6} and \eqref{xi7} imply
\begin{align}
  \label{xi8}
  \J(\uhk) =\ell_h(\zhk),
\end{align}
which together with \eqref{eq:eqi} imply \eqref{eq:equi_algB}. We recall that the last
relation is valid only if $\uhk$ and $\zhk$ are obtained by the BiCG method in exact arithmetic.

\subsection{Estimation of the algebraic error}
\label{sec:estim}

In Sections~\ref{sec:eqi} and \ref{sec:bicg}, we presented  three
  possibilities of the evaluations of the quantity of interest {\PPa} -- {\PPc}.
  In this paragraph, we discuss the estimates of the errors of these evaluations,
  i.e., estimation of quantities $\estAk$ and $\estDAk$ introduced
    in \eqref{JP2}  and in \eqref{JP5}, respectively.    
We use the standard approach when $\nu>0$
additional algebraic solver steps are performed and the difference
between $k$ and $k+\nu$ iterates is used for the estimate of the error.

%
%
\begin{itemize}
\item [\EEa] To estimate the error $\bcT \mA^{-1} \brk$ of evaluation  {\PPa}, we subtract the identity \eqref{BB1} in iterations $k$ and $k+\nu$, and obtain
  \begin{align}
  \label{EE3}
    \bcT \mA^{-1} \brk
    =  \underbrace{\bcT(\bxknu-\bxk)}_{\mathrm{estimate}}
    + \,\bcT \mA^{-1} \brknu.
  \end{align}


  Similarly, the dual identity \eqref{BB2} gives
  \begin{align}
    \bskT \mA^{-1} \bb
    =  \underbrace{(\byknu-\byk)^\T \bb}_{\mathrm{estimate}}
    + \,\bsknuT \mA^{-1} \bb.
  \end{align}  
\item [\EEb] Subtracting the identities \eqref{BB3} in iterations $k$ and $k+\nu$ we can express the error $\bskT \mA^{-1} \brk$ of evaluation {\PPb} as
  \begin{align}
     \label{EE6}
    \bskT \mA^{-1} \brk
    =  \underbrace{\bcT(\bxknu-\bxk)+\byknuT\brknu-\bykT\brk}_{\mathrm{estimate}}
    + \,\bsknuT \mA^{-1} \brknu.
  \end{align}
  Analogously, using \eqref{BB4} we obtain
  \begin{align}
    \bskT \mA^{-1} \brk
    =  \underbrace{(\byknu-\byk)^\T \bb+\bsknuT\bxknu-\bskT\bxk}_{\mathrm{estimate}}
    + \,{\bsknuT \mA^{-1} \brknu}.
  \end{align}  

\item [\EEc] 
  Finally, considering the identity \eqref{xi5} in iterations $k$ and $k+\nu$ we can express the error $\bskT \mA^{-1} \brk$ of evaluation {\PPc} as 
  \begin{align}
    \label{EE16}
     \bskT \mA^{-1} \brk
   = \underbrace{\bxiBknu - \bxiBk }_{\mathrm{estimate}}
    +\, \bsknuT \mA^{-1} \brknu.
  \end{align}
\end{itemize}
Similarly as in Section~\ref{sec:eqi}, starting with $\bx_0=0$ and $\by_0=0$, 
all errors of evaluations {\PPa}--{\PPc} as well as their estimates {\EEa}--{\EEc} are identical for the BiCG method 
in exact arithmetic. 
However, in finite precision arithmetic when the
orthogonality \eqref{orthog} is violated, 
the error of the evaluation {\PPa} and its estimate {\EEa}
can substantially differ from the errors of the evaluations  {\PPb}--{\PPc}
and their estimates {\EEb}--{\EEc}.

\section{Numerical experiments on fixed meshes}
\label{sec:numerF}

In this section we present the first collections of numerical
experiments where approximation space $\Vh$ is fixed.
The aim is to demonstrate the accuracy of the approximation
of the quantity of interest $\J(\uh) = \bcT \mA^{-1} \bb = \ell_h(\zh)$
from Sections~\ref{sec:eqi} and \ref{sec:bicg}, and the estimates of the errors of these
these approximations from Section~\ref{sec:estim}.

We consider a second order elliptic problem which is discretized by 
the symmetric interior penalty Galerkin (SIPG) method using
a piecewise polynomial but discontinuous approximation.
SIPG method guarantees the
the primal as well as dual consistencies \eqref{eq:cons} and \eqref{eq:consD},
respectively. For the definitions of the forms $\ah$ and $\ell_h$,
we refer to \cite{ESCO-18}, the detailed analysis can be found, e.g., in \cite{DGM-book}.
All numerical examples presented in this paper
were carried out using the in-house code ADGFEM \cite{ADGFEM} written
in gfortran in double precision with processor i7-2620M CPU 2.70GHz (Ubuntu 16.04).

\subsection{Elliptic problem on a ``cross'' domain}
\label{sec:cross}

\begin{figure}[t]
  \includegraphics[height=0.33\textwidth]{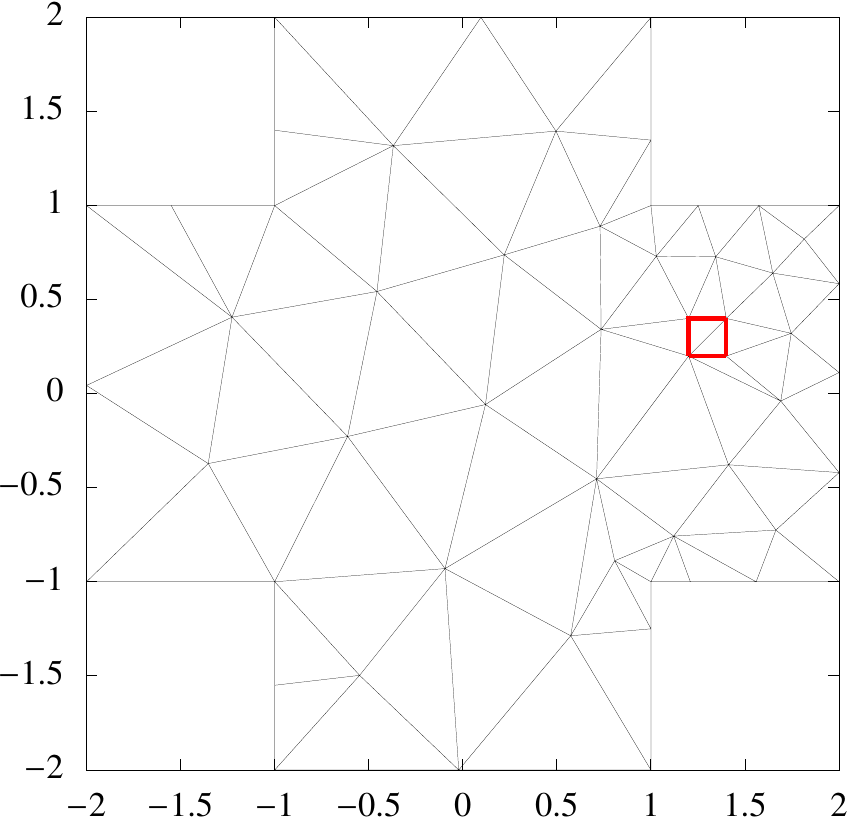}
  \includegraphics[height=0.33\textwidth]{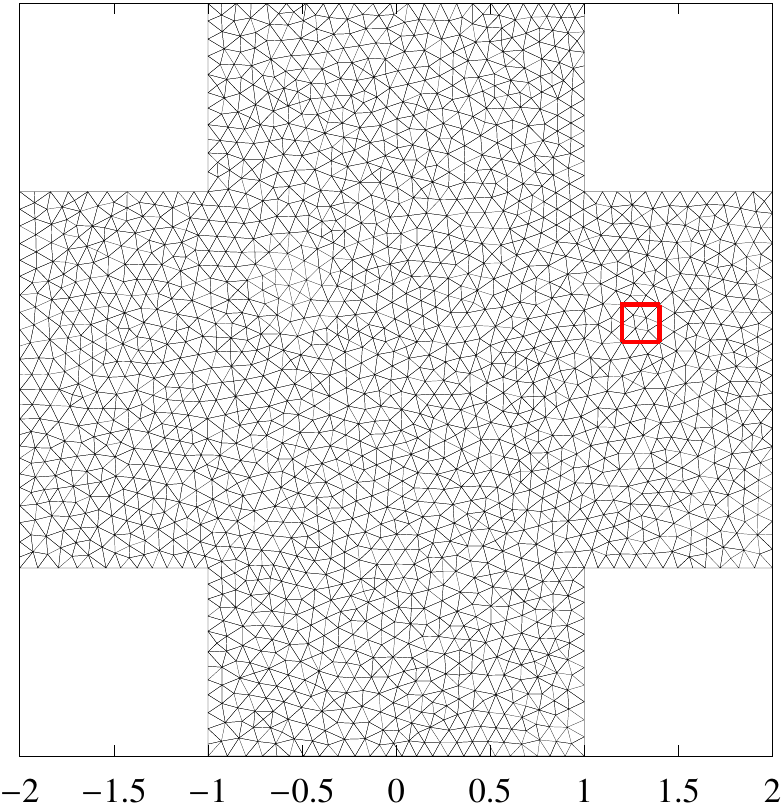}
  \includegraphics[height=0.33\textwidth]{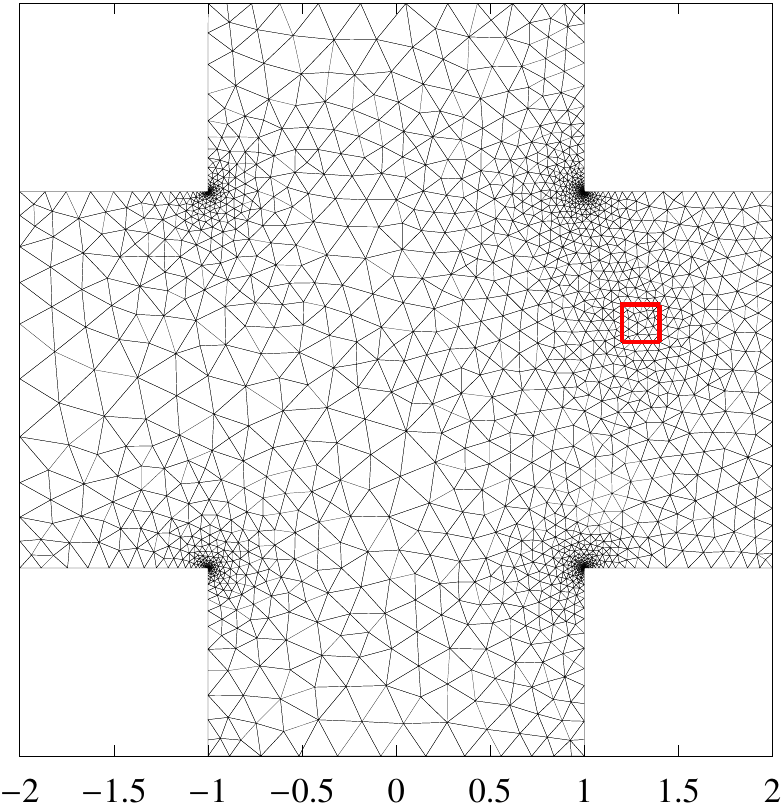}
  \caption{Cross domain, computational domain $\Om$ an initial mesh and
    the domain of interest $\Om_J$ (small red square)  (left),
    the finer uniform mesh (center) and the adaptively refined mesh (right).}
  \label{fig:CR_meshes}
\end{figure}

We consider the example from  \cite[Example 2]{Ainsworth2012Guaranteed}
\begin{align}
  \label{eq:CR}
 - \Delta u  = 1 \quad \text{ in } \Om,  \qquad 
  u  = 0  \quad \text{ on } \partial\Om,
\end{align}
where the ``cross'' domain $\Om = (-2,2) \times (-1,1) \cup (-1,1) \times (-2,2)$
and  $\Delta$ denotes the Laplace operator. 
The target functional is defined as the mean value of the solution over the square
$\Om_J = [1.2,1.4]\times [0.2,0.4],$ i.e. 
$J(u) = \frac{1}{|\Om_J|}\intX{\Om}{}{ j_\Om(x) u(x) }{x},$
where $j_\Om$ is the characteristic function of the square $\Om_J$,
see Figure~\ref{fig:CR_meshes}, left.
The exact value of $J(u)$ is unknown but we use the reference value $0.407617863684$,
which was computed in \cite{Ainsworth2012Guaranteed} on an
adaptively refined mesh with more than 15 million triangles. 

The presence of interior obtuse angles of $\Om$ gives the singularities
of the weak solution of \eqref{eq:CR}.
We carried out the computations on two triangular meshes,
the first one is (quasi-)uniform having 3742 triangles
and the second one is adaptively refined in the vicinity of interior angles
and it has 4000 triangles, see Figure~\ref{fig:CR_meshes}.
For both meshes we used the SIPG method with $P_2$ and $P_4$ polynomial approximations

\subsection{Approximation of the quantity of interest}
\label{sec:crossR}

For each of the four corresponding discrete problems (uniform/adapted mesh and
$P_2$/$P_4$ approximations) we carried out the solution of the corresponding
algebraic systems \eqref{AA1} by the BiCG method from Algorithm~\ref{alg:BiCG}.
Since the method has tendency to stagnate after some number of iterations,
we restarted the computations once after 400 BiCG iterations. 
Table~\ref{tab:CR} shows the limit values of
$\J(\uh) = \bcT\mA^{-1}\bb$.

\begin{table}
  \begin{center}
    \bgroup
    \def\arraystretch{1.25}
    \begin{tabular}{cc|c|c}
      \hline
      mesh    uniform & $P_2  $ & $\DoF=$     22452 & 0.4071783143507507 \\
      mesh    uniform & $P_4  $ & $\DoF=$     56130 & 0.4075262691478035 \\
      mesh    adapted & $P_2  $ & $\DoF=$     24000 & 0.4076152998044911  \\
      mesh    adapted & $P_4  $ & $\DoF=$     60000 & 0.4076169203362077  \\
      \hline
    \end{tabular}
    \egroup
  \end{center}
  \caption{Cross domain, computation on fixed meshes,
    the limit values of the quanitity of interest $\J(\uh) = \bcT\mA^{-1}\bb$
    after three restarts on two meshes (uniform/adapted)
    and $P_2$ and $P_4$ approximations.}
  \label{tab:CR}  
\end{table}

\begin{figure}[t]
  \includegraphics[width=1.0\textwidth]{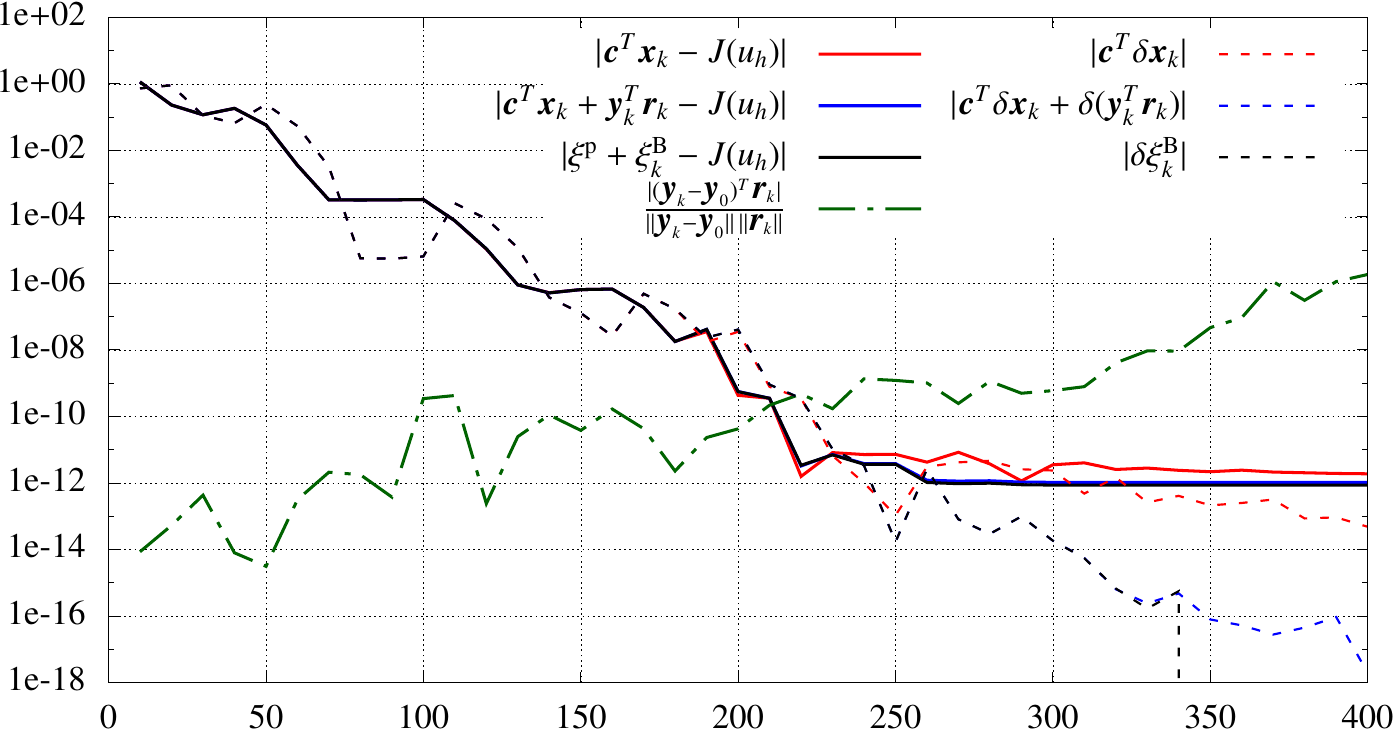}
   \caption{Convergence of the algebraic error and its estimates for the BiCG solver,
    case uniform mesh, $P_2$ approximation,
    the errors of the approximation of the quantity of interest {\PPa}, {\PPb}, {\PPc}
    (solid lines),
    their approximation {\EEa}, {\EEb}, {\EEc} (dashed lines) and the
    lost of the orthogonality (dotted-dashed line).}
  \label{fig:BiCG_convU1P2}
\end{figure}

\begin{figure}[t]
  \includegraphics[width=1.0\textwidth]{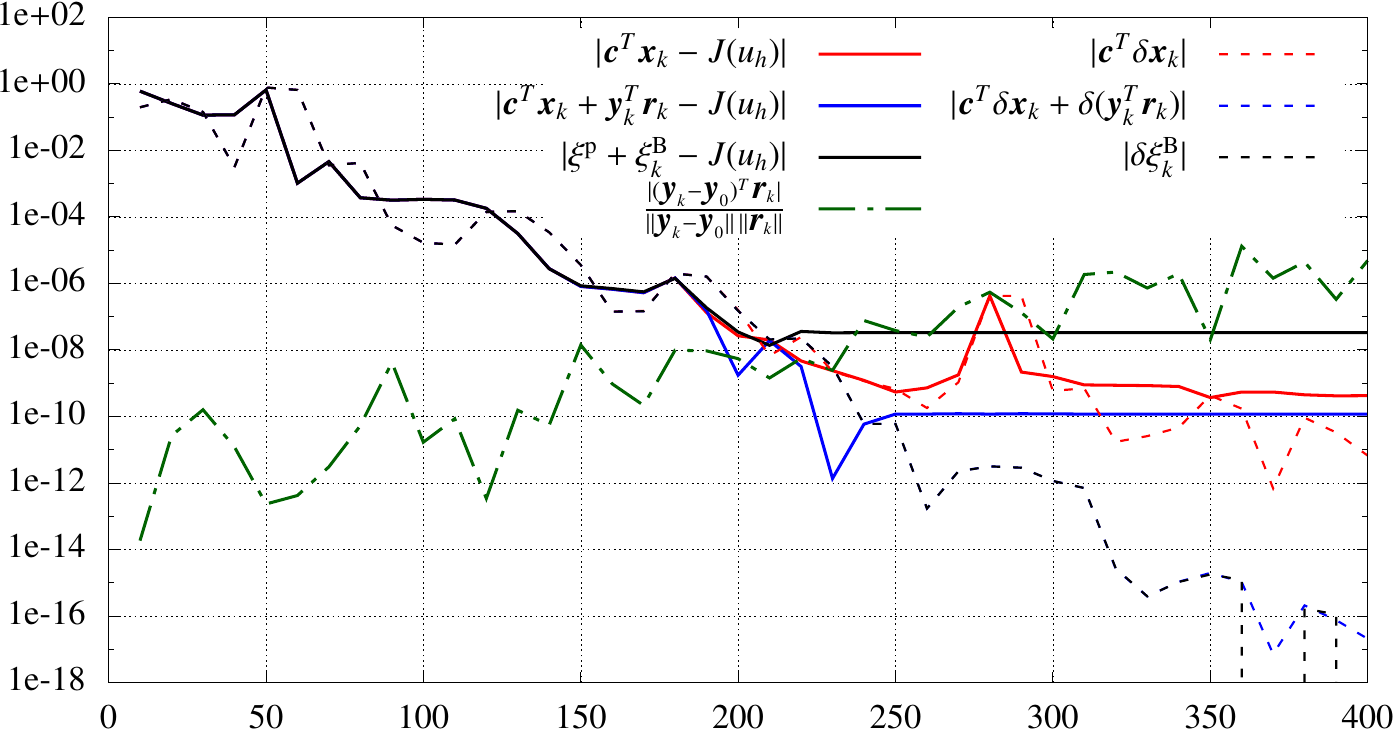}
  \caption{Convergence of the algebraic error and its estimates for the BiCG solver,
    case uniform mesh, $P_4$ approximation,
    the errors of the approximation of the quantity of interest {\PPa}, {\PPb}, {\PPc}
    (solid lines),
    their approximation {\EEa}, {\EEb}, {\EEc} (dashed lines) and the
    lost of the orthogonality (dotted-dashed line).}
  \label{fig:BiCG_convU1P4}
\end{figure}

\begin{figure}[t]
  \includegraphics[width=1.0\textwidth]{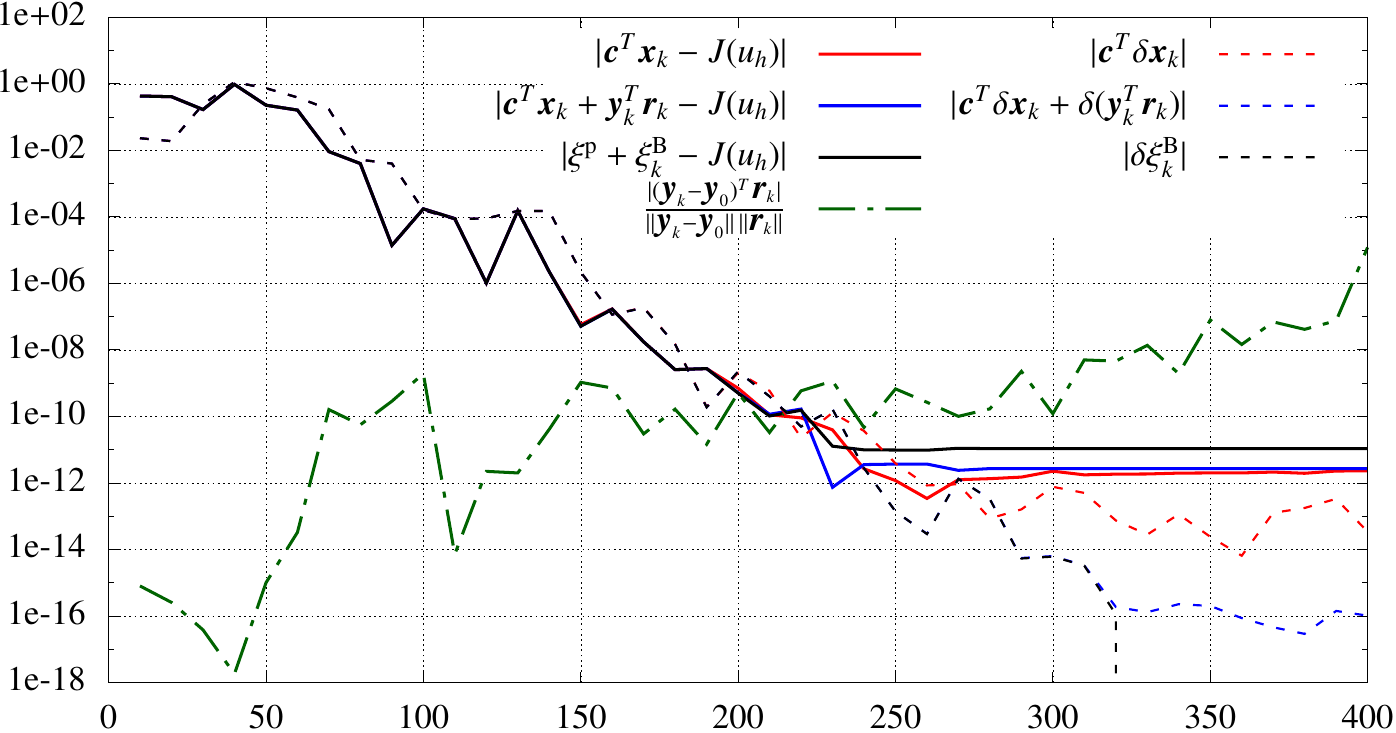}
  \caption{Convergence of the algebraic error and its estimates for the BiCG solver,
    case adapted mesh, $P_2$ approximation, 
    the errors of the approximation of the quantity of interest {\PPa}, {\PPb}, {\PPc}
    (solid lines),
    their approximation {\EEa}, {\EEb}, {\EEc} (dashed lines) and the
    lost of the orthogonality (dotted-dashed line).}
  \label{fig:BiCG_convA2P2}
\end{figure}

\begin{figure}[t]
  \includegraphics[width=1.0\textwidth]{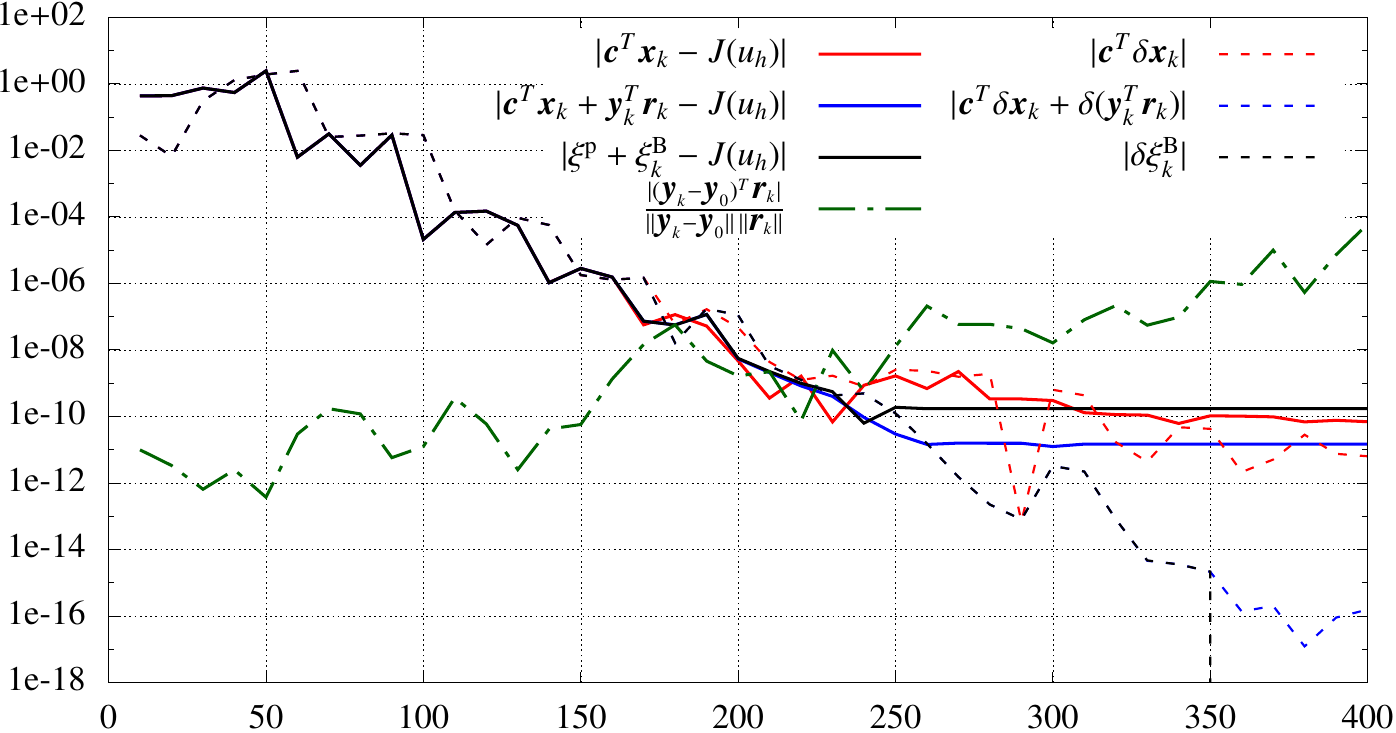}
  \caption{Convergence of the algebraic error and its estimates for the BiCG solver,
    case adapted mesh, $P_4$ approximation, 
    the errors of the approximation of the quantity of interest {\PPa}, {\PPb}, {\PPc}
    (solid lines),
    their approximation {\EEa}, {\EEb}, {\EEc} (dashed lines) and the
    lost of the orthogonality (dotted-dashed line).}
  \label{fig:BiCG_convA2P4}
\end{figure}

Figures~\ref{fig:BiCG_convU1P2} -- \ref{fig:BiCG_convA2P4} show
  the results obtained within the first 400 BiCG
  iterations\footnote{After the restart, the machine accuracy is achieved in few steps and the error estimators give vanishing values, therefore we do not show them.}, namely
\begin{itemize}
\item[$\bullet$] the convergence
  of the errors  three types of the approximation of the
  quantity of interest $\J(\uh) = \bcT\mA^{-1}\bb$, namely
  \begin{itemize}
    \setlength\itemsep{2pt}
  \item {\PPa}: $ | \bcT\bxk - \J(\uh) | $ following from \eqref{BB1},
  \item {\PPb}: $ | \bcT\bxk + \bykT \brk -  \J(\uh) | $ following from \eqref{BB3},
  \item {\PPc}: $ | \bxiP + \bxi_{k}^{B} -  \J(\uh) | $ following from \eqref{xi5},
  \end{itemize}
  where $\J(\uh)$ is the limit value from Table~\ref{tab:CR} and $|\cdot|$ is the absolute value,
\item[$\bullet$]  the convergence of the values of three different estimations
  of the error of the approximation of the quanitity of interest, namely
  \begin{itemize}
    \setlength\itemsep{2pt}
  \item {\EEa}: $ | \bcT\delta \bxk |$ following from \eqref{EE3},
  \item {\EEb}: $ | \bcT \delta \bxk +\delta(\bykT\brk) |$ following from \eqref{EE6},
  \item {\EEc}: $ | \delta \bxiBk|$ following from \eqref{EE16},
  \end{itemize}
  where symbol $\delta$  denotes the difference between $k$-th and $(k+\nu)$-th iterations,
  i.e,  $\delta \bxk= \bxknu - \bxk$,
  $\delta(\bykT\brk)= \byknuT\brknu-\bykT\brk$, 
  $\delta \bxiBk = \bxiBknu - \bxiBk$, in the experiments, we put uniquely $\nu=10$,
\item[$\bullet$]  the quantity
  \begin{itemize}
    \setlength\itemsep{2pt}
  \item $ \frac{| (\byk-\by_0)^\T \brk  |}{\|\byk-\by_0 \|\, \| \brk\|} $
    measuring the lost of the orthogonality, cf.  \eqref{orthog}.
  \end{itemize}
\end{itemize}

{
We observe the following. 
\begin{enumerate}
\item The errors of all approximations of the quantity of interest (solid lines) 
  decrease until they reach some level of accuracy. After reaching this level they stagnate.
  This is in agreement with the results in 
  \cite{Gr1997} where its is showen that the residual norms $\| \brk \|$
  of the Krylov subspace methods like BiCG, reach (if they converge)
  the level of accuracy close to $\varepsilon \|A\|  \Theta $
  where $\varepsilon$ is the machine precision and $\Theta = \max\{\|\bxk\|\}$ 
  (unpreconditioned case). 
\item All approximation of the quantity of interest {\PPa} -- {\PPc} (solid lines)
  have almost identical convergence as long as the lost of orthogonality \eqref{orthog}
  is small (dotted-dashed line). When  the lost of orthogonality starts to play
  more important role the values of the approximation of the quantity of interest
  slightly differ.
\item Similarly, all the estimations of the error of approximations
  {\EEa} -- {\EEc} (dashed lines)
  have almost identical convergence as long as the lost of orthogonality
  is small. If the evaluations of the quantity of interests start to stagnate
  the estimates  {\EEa} -- {\EEc} underestimate the error. This process is the slowest
  for {\EEa} but it is presented too.
\end{enumerate}
From these observations, we conclude that any error estimate {\EEa} -- {\EEc} can
be used for the stopping criterion. They underestimate the error only when the
computational process stagnates and then it make no sense to proceed with
next iterative steps. Based on the argumentation in Section~\ref{sec:estim},
criteria {\EEb} and {\EEc} are less sensitive to the lost of orthogonality.
Finally, criterion {\EEc} is cheaper to evaluate than {\EEb} but, in our case,
it does not play
any essential role in the whole computational process.
}

\section{Adaptive mesh refinement and algebraic stopping criteria}
\label{sec:adapt}

The goal of the numerical solution of \eqref{eq:PP} is to obtain a numerical
approximation $\tuh$ such that (cf. \eqref{EEP})
\begin{align}
  \label{goal}
  |\J(\u) - \J(\tuh)| \approx   \etaSk + \etaAk \leq \TOL,
\end{align}
where $\TOL>0$ is the given tolerance. The adaptive mesh refinement allows to
reduce the computation  costs necessary to achieve \eqref{goal}.

\subsection{Mesh adaptation algorithm}

The idea of the mesh adaptive algorithm is to start on an initial coarse
mesh $\Thz$ (dimension of the corresponding space $\Vhz$ is small).
Then for $m=0,1,\dots$, we discretize and solve both primal and dual problem
on $\Vhm$ and estimate the error of the quantity of interest.
If the estimate does not fulfil \eqref{goal} then we adapt the mesh and
create a new one $\ThmP$ and proceed with the computation.
Algorithm~\ref{alg:AMAhp} shows the abstract form of the adaptive mesh refinement algorithm.
It means that we obtain the approximations of the primal and dual
solutions  $\uhmk\in\Vhm$ and $\zhmk\in\Vhm$, where
subscript $m$-th corresponds to the level of mesh adaptation and
subscript $k$-th corresponds to the algebraic iteration.


\begin{algorithm}[ht]
  \caption{Mesh adaptation process}
  \label{alg:AMAhp}
  \begin{algorithmic}[1]
{
  \STATE let $\Thz$ be the initial (coarse) mesh
  \STATE let $\Vhz$ be the
  corresponding finite element space 
\FOR{$m=0,1,\dots$} 
  \STATE set the algebraic problems \eqref{AA1} corresponding to the
  discretization of \eqref{eq:PP} and \eqref{eq:DP}, respectively, on $\Vhm$, i.e.,
  \begin{align}
  \label{AAm}
  \mA_m \bx_{m} = \bb_{m}\qquad\mbox{and}\qquad \mA_{m}^{\T}\by_m=\bc_m,
\end{align}
  where $\bb_{m},\bc_m \in\R^{N_m}$, $\mA_m\in\R^{N_m\times N_m}$,  $\dim\Vhm=N_m$,
  \STATE apply an iterative algebraic solver 
  \FOR{$k=0,1,\dots$}
  \STATE evaluate $\bx_{m,k}$ and $\by_{m,k}$ approximations of \eqref{AAm}
  \IF{ {\em algebraic stopping criterion} is achieved}  \STATE EXIT \ENDIF
  \ENDFOR
  \STATE define $\uhmk,\zhmk\in\Vhm$ corresponding
  to $\bx_{m,k}$ and $\by_{m,k}$ (the outputs of Algorithm~\ref{alg:BiCG})
  \STATE define  higher-order reconstructions
  $\huhm=\krr(\uhmk)$ and $\hzhm=\krr(\zhmk)$ using \eqref{krr}
\STATE employing \eqref{etas}--\eqref{etasA}, evaluate the goal oriented error estimates
  \begin{align}
    \label{etasm}
    \etamk := \etaSmk + \etaAmk, \qquad \etaDmk := \etaDSmk + \etaDAmk
  \end{align}
  where 
  \begin{align}
    \label{etasm1}
    &\etaSmk := \res{\uhmk}{\hzhm -  \zhmk},  &&\etaAmk := \res{\uhmk}{\zhmk}, \\
    &\etaDSmk := \resD{\zhmk}{\huhm -  \uhmk},  &&\etaDAmk := \resD{\zhmk}{\uhmk}, \notag
  \end{align}
\IF {$\etamk \leq \TOL$ and $\etaDmk \leq \TOL$} \STATE{ STOP the computational process 
  with the output quantity $\J(\uhmk)$ and its estimate $\etamk$},
\ELSE \STATE 
  based on a localization of $\etamk$, adapt the mesh, i.e., modify the size and/or
  the shape of elements and possibly also the polynomial approximation degrees,
  the new mesh is $\ThmP$ and the corresponding space $\VhmP$,
  \ENDIF
  \ENDFOR
  }
\end{algorithmic}
\end{algorithm}

Step 18 of Algorithm~\ref{alg:AMAhp} (mesh adaptation)
depends on the used discretization method and the refinement technique.
In order to demonstrate the robustness of the presented stopping criteria,
we use the anisotropic $hp$-mesh adaptation approach from \cite{ESCO-18},
which generate anisotropic meshes (consisting of possibly thin and long
triangular elements) and varying polynomial approximation degree.

The crucial aspect is the algebraic stopping criterion in step 8
of Algorithm~\ref{alg:AMAhp} which is discussed in the next section.

\subsection{Standard stopping criteria for the solution of algebraic systems}
\label{sec:stop1}
We focus on the algebraic stopping criterion in Step 8 of Algorithm~\ref{alg:AMAhp}.
Obviously, too strong criterion leads to many algebraic iterations without a gain
of accuracy. On the other hand, the weak criterion leads to an under-solving
of \eqref{AAm} which affect the mesh adaptation process. Typically too many
mesh elements are generated.

Often the residual stopping criteria for \eqref{AAm} are used, i.e.,
\begin{align}
  \label{RES}
  \| \br_{m,k} \| := \| \bb_m - \mA_m \bx_{m,k} \| \leq \TOLA, \qquad
  \| \bs_{m,k} \| := \| \bc_m - \mA_m^{\T} \by_{m,k} \| \leq \TOLA, 
\end{align}
or their preconditioned variant
\begin{align}
  \label{RESP}
  \| \mP^{-1} (\bb_m - \mA_m \bx_{m,k}) \| \leq \TOLA, \qquad
  \| \mP^{-\T} (\bc_m - \mA_m^{\T} \by_{m,k}) \| \leq \TOLA, 
\end{align}
where $\mP$ represents a suitable preconditioner, cf.~\eqref{mP}.
This criterion is easy to evaluate but the choice of the suitable tolerance $\TOLA>0$
is difficult since this criterion has no relation to the discretization error.

This drawback was eliminated in \cite{Meidner2009Goal} by the following
stopping criterion
\begin{align}
  \label{aDWR}
  \etaAmk \leq \cA \etaSmk,\qquad
  \etaDAmk \leq \cA \etaDSmk,
\end{align}
where the primal/dual estimates of the algebraic/discretization errors are defined
by \eqref{etasm1} and $\cA\in(0,1)$ is a suitable constant.
This means that Steps 12-14  of Algorithm~\ref{alg:AMAhp} are moved inside
the inner loop (after Step 7).
We call this criterion as {\em algebraic goal-oriented} stopping criterion.

The conditions \eqref{aDWR} allow to control the size of the algebraic error.
However, a strong drawback of \eqref{aDWR} are to computational costs.
Whereas the evaluation of $ \etaAmk $ and $ \etaDAmk$ is cheap,
the computation of $ \etaSmk $ and $ \etaDSmk$ is much more expensive namely due to
the necessity to perform the higher-order reconstructions $\huhm$ and $\hzhm$
(Step 13).
The computational costs can be reduced by testing \eqref{aDWR}
only after some number of iterations, e.g., after 1-3 restarts of iterative solver.
In \cite{Meidner2009Goal}, the first condition of \eqref{aDWR} was tested
after the performing of one cycle of multigrid method.

\subsection{New stopping criteria for the solution of algebraic systems}
\label{sec:stop2}

In Section~\ref{sec:alter}, we
  introduced two alternative formula for the decomposition of the computational errors
  into the discretization and algebraic parts, namely using \eqref{JP1} and \eqref{JP2},
  we have
  \begin{alignat}{3}
    \label{ST1}
    \J(\u - \uhmk) &  = \res{\uhmk}{\z -  \zhmk} + \res{\uhmk}{\zhmk} &&=\errSmk +\errAmk,
    \\
    \J(\u - \uhmk) &  =  \J(\u - \uhm) + \J(\uhm - \uhmk)&&= \estSmk +\estAmk. \notag
  \end{alignat}
  Similarly, \eqref{eq:JD_alg} and \eqref{JP5} imply the dual counterpart
  \begin{alignat}{3}
    \label{ST2}
    \ell_h(\z - \zhmk) & = \resD{\zhmk}{\u -  \uhmk} + \resD{\zhmk}{\uhmk} &&= \errDSmk + \errDAmk, \\
    \ell_h(\z - \zhmk) &= \ell_h(\z - \zhm) + \ell_h(\zhm - \zhmk)   &&= \estDSmk + \estDAmk.
    \notag
  \end{alignat}
  The discussion presented therein shows that both quantities $\errAmk\not=\estAmk$
  corresponds to the algebraic error of the primal  problem and similarly,
  $\errDAmk\not=\estDAmk$
  corresponds to the algebraic error of the dual one.
  
In Section~\ref{sec:estim}, we presented
several techniques estimating the quantities $\estAmk$ and $\estDAmk$ by techniques
{\EEa} -- {\EEc}. On the other hand, quantities $\errAmk$ and  $\errDAmk$ are
available during the BiCG iterative method since
\begin{align}
  \label{ST3}
  \errAmk &= \etaAmk = \res{\uhmk}{\zhmk} = \bykT \brk, \\
  \errDAmk &= \etaDAmk = \resD{\zhmk}{\uhmk} = \bskT \bxk, \notag
\end{align}
cf. \eqref{etasA} and \eqref{aa5}.

\begin{figure}[t]
  \begin{center}
    \includegraphics[width=\textwidth]{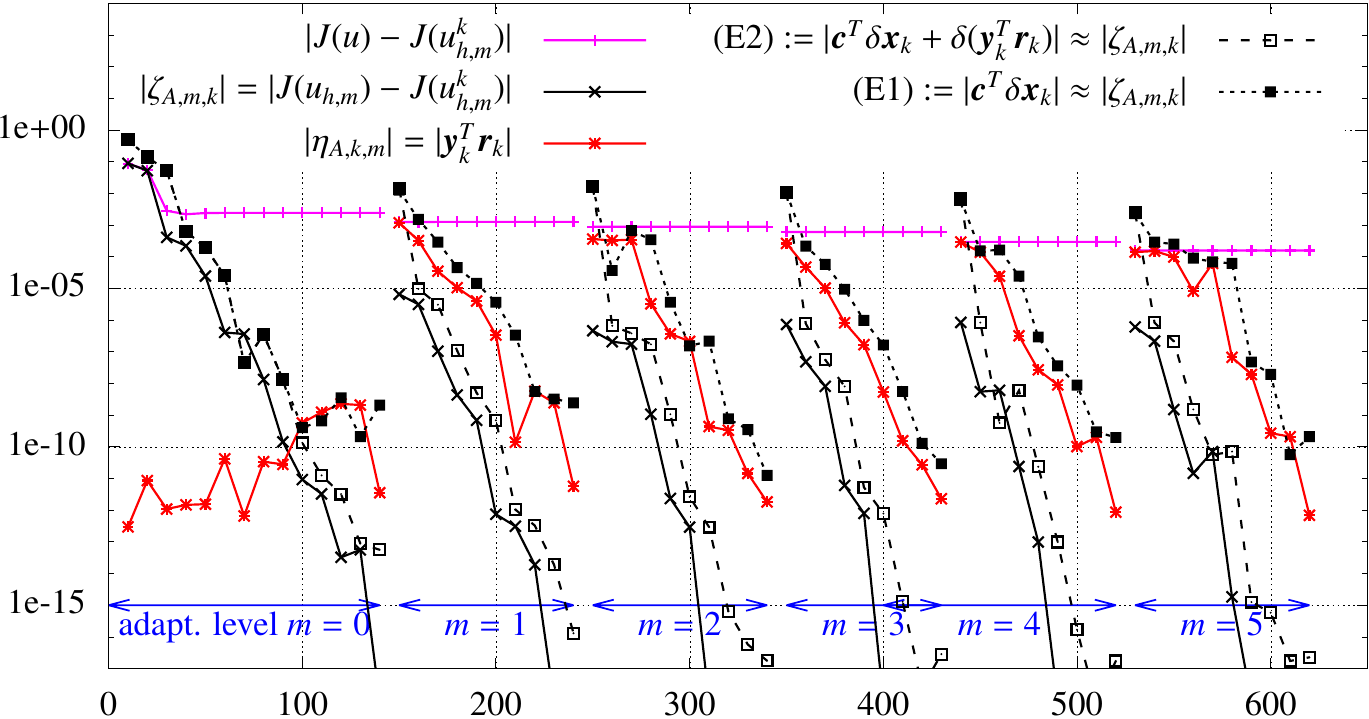}

    \vspace{4mm}
    
    \includegraphics[width=\textwidth]{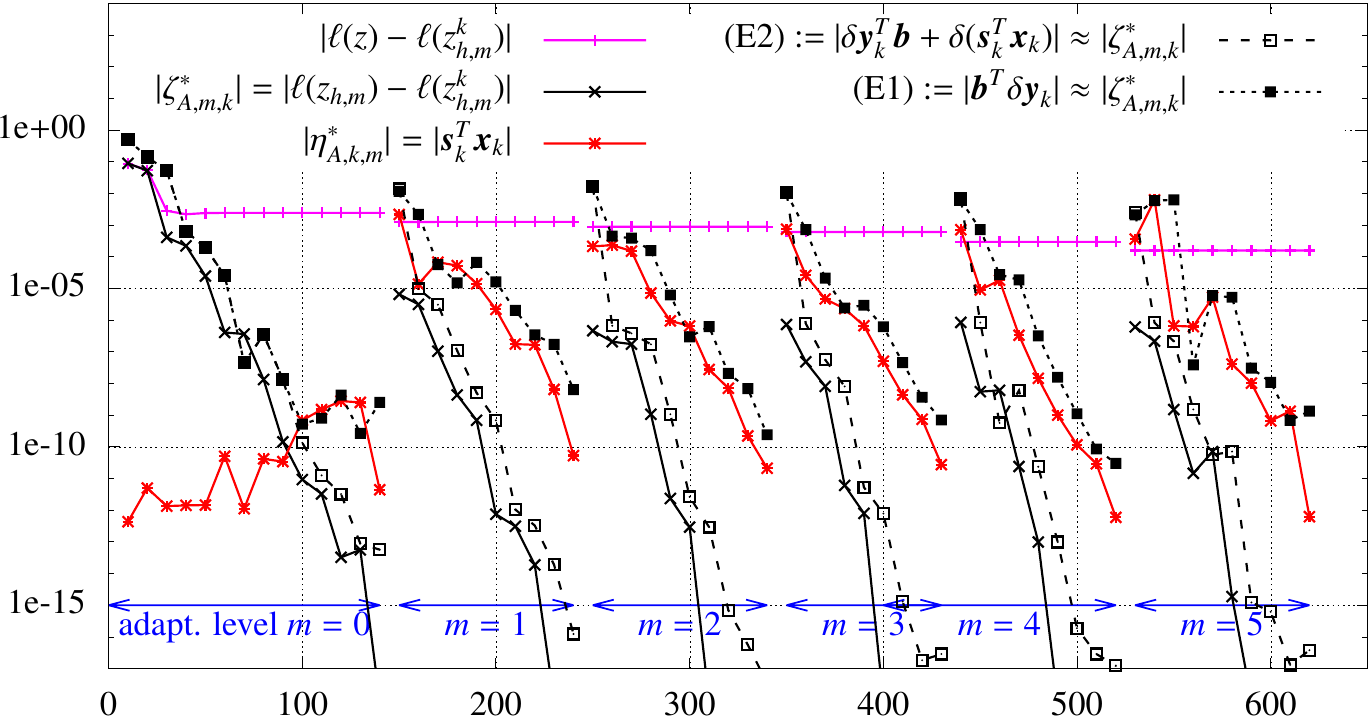}
  \end{center}
  \caption{Typical example of the adaptive computations for mesh adaptation
    steps $m=0,\dots,5$,
    dependence of primal (top) and dual (bottom) quantities
    w.r.t. the number of BiCG steps.}
  \label{fig:adaptBiCG}
\end{figure}

Figure~\ref{fig:adaptBiCG} shows a typical dependence of these quantities
w.r.t. the number of BiCG iterations
during mesh adaptation process by  Algorithm~\ref{alg:AMAhp} for $m=0,\dots, 5$.
In the top figure,
we plot the total error $J(\u-\uhmk)$, the exact algebraic error $\estAmk=J(\uhm-\uhmk)$,
the algebraic error $\errAmk=\bskT \bxk$ and the estimates {\EEa} and {\EEb}
of $\estAmk$ from Section~\ref{sec:estim}.
The bottom figure shows the dual counterparts. 
Estimate {\EEc} gives the same graphs as {\EEb} so we do not show it.

We observe that estimate  {\EEb} approximates $\estAmk$
similarly as {\EEa} for the moderate values of accuracy but
much better on the level close to the machine accuracy.
This
is in agreement with theoretical considerations in Section~\ref{sec:estim},
but this effect is not observed on the initial mesh ($m=0$)
when we do not have a good initial approximation.
Further, both algebraic error representations
converge but with different speed. Whereas on the initial mesh we have
$\eta_{\mathrm{A},0,k} \ll  \zeta_{\mathrm{A},0,k}$ for $k\le 100$,
starting from $m=1$ we observe $\eta_{\mathrm{A},m,k} \gg  \zeta_{\mathrm{A},m,k}$.
Similar behaviour is observed for the dual quantities.



Based on this observation, we consider both pairs of quantities
$(\etaAmk,  \estAmk)$ and $(\etaDAmk,  \estDAmk)$ for the definition of the
stopping criterion, namely,  we evaluate the quantities
\begin{align}
  \label{estMK1}
  \ESTAmk := & | \delta \bxiBk| + | \etaAmk| = |\bxiBknu - \bxiBk| + | \by_{k,m}^{\T} \br_{k,m}|, \\
  \ESTDAmk := & | \delta \bxiBk| + | \etaDAmk| = |\bxiBknu - \bxiBk| + | \bs_{k,m}^{\T} \bx_{k,m}|, \notag
\end{align}
corresponding to the $k$-iteration  of the BiCG Algorithm~\ref{alg:BiCG}
on $m$-level of mesh adaptation.
Let us recall that quantities $\ESTAmk$ and $\ESTDAmk$ are available during
the whole BiCG iterative process with negligible computational costs.

Now we define a new algebraic stopping criterion
(called hereafter $\sigma$-stopping criterion)
for Algorithm~\ref{alg:AMAhp}
(Step~8) by
\begin{align}
  \label{xiB}
  \ESTAmk  \leq \cA \TOL\qquad \mbox{and} \qquad
  \ESTDAmk \leq \cA \TOL, 
\end{align}
where $ \TOL>0$ is the (global) tolerance from \eqref{goal}
and $\cA\in(0,1)$ is a suitable constant.
In contrary to the stopping criterion \eqref{aDWR}, the proposed
one \eqref{xiB} does not take into account the (estimate of) the discretization
error and consequently the strong advantage of \eqref{xiB} is that
its evaluation is very fast.

We can summarized the idea of the original stopping criterion \eqref{aDWR} as follows:
``{\it algebraic system is solved as long the algebraic error is $\cA$-times smaller
than the discretization error}''.
On the other hand, the idea of the new stopping criterion \eqref{xiB} is the following: 
``{\it algebraic system is solved as long the algebraic error is $\cA$-times smaller
than the prescribed tolerance for the discretization error}''.

The stopping criterion  \eqref{xiB} does not look to much efficient since for
starting levels of mesh adaptation, when $\etaSmk \gg\TOL$ and $\etaDSmk \gg\TOL$,
the iterative solver is stopped when
the algebraic error is much lower than the discretization one.
It is true but on the other hand performing of additional several teens of BiCG iterations
is typically faster than the evaluation of $\etaSmk$ and $\etaDSmk$.
Moreover, approximate solutions computed using BiCG on starting levels serve as good initial approximations for BiCG on higher levels, which improve then the solving process significantly.


\section{Numerical experiments on adaptively refined meshes}
\label{sec:numerA}

We demonstrate the computational performance of the stopping criteria
from Sections~\ref{sec:stop1} and \ref{sec:stop2}.
We consider two numerical examples, the first one is the elliptic problem
on the ``cross'' domain described in Section~\ref{sec:cross} and
the second one is a convection-dominated problem having some anisotropic features,
it is defined in Section~\ref{sec:carpio}.
Both problems are discretized again by the SIPG method
and the meshes are adapted by the technique from \cite{ESCO-18}.
Each mesh $\Thm$ is generated from the computed primal and dual solutions on the
previous one. These solutions are interpolated to the actual mesh, hence
we have relatively very good initial approximations $\bx_{m,0}$ and $\by_{m,0}$
for the solution of \eqref{AAm}.

\subsection{Convection-dominated problem}
\label{sec:carpio}

The second example is 
taken from \cite{FormaggiaPerotto04}, see also \cite{ESCO-18,CarpioPrietoBermejo_SISC13}.
We solve the convection-diffusion equation
\begin{align}
  \label{carpio}
  - \ve \Delta u + \nabla \cdot (\bkb\, u) &= 0  \qquad
  \text{in } \Om:=[0, 4]\times [0, 4] \setminus [0, 2] \times [0, 2],
\end{align}
where $\ve = 10^{-3}$,
the convection field $\bkb = (x_2,-x_1)$
and $\nabla \cdot$ is the divergence operator.
We prescribe the Dirichlet and Neumann boundary conditions
\begin{align}
  \label{carpio:BC}
  u=1 &\quad  \mbox{ on } \{x_1 = 0\}, \\
  \nabla u \cdot \bkn = 0 & \quad \mbox{ on } \Gamma_1:=\{\gom; x_1 = 4\} \cup \Gamma_2:=\{\gom; x_2 = 0\},  \notag \\
  u=0 & \quad \mbox{ elsewhere}. \notag
\end{align}
The solution $u$ exhibits boundary layers as well as two circular-shaped internal
layers.
We  consider the functional 
$  J(u) = \int_{\Gamma_1} \bkb\cdot\bkn\, u\dS$
with the reference value $J(u) = 0.07408122\pm 10^{-8}$.

\subsection{Comparison of the stopping criteria}

In order to demonstrate the robustness of the proposed stopping criteria,
we solve the elliptic problem \eqref{eq:CR} and
the convection-diffusion problem \eqref{carpio}  
by the mesh adaptation Algorithm~\ref{alg:AMAhp}.
For both problems, we employ the anisotropic $h$-mesh adaptation using $P_3$ approximation
and anisotropic $hp$-mesh adaptation, for details we refer to \cite{ESCO-18}.
In order to observe the effect with not sufficiently resolved algebraic systems
\eqref{AAm}, we stop the adaptation process when 
\begin{align}
  \label{STOP}
  \frac12\left(\etaSmk + \etaDSmk\right)\leq \TOL,
\end{align}
where $\etaSmk$ and $\etaDSmk$ are given by \eqref{etasm1} and we
put $\TOL=10^{-8}$ for the elliptic problem \eqref{eq:CR}
and $\TOL=10^{-10}$ for the convection-diffusion problem \eqref{carpio}.

We test Algorithm~\ref{alg:AMAhp} with the following iterative solvers and
the stopping criteria.
\begin{enumerate}[label=(\alph*)]
\item GMRES with preconditioned residual stopping criterion \eqref{RESP}
  with tolerances $\TOLA=10^{-9}$, $10^{-6}$ and $10^{-3}$,
\item BiCG with preconditioned residual stopping criterion \eqref{RESP}
  with tolerances $\TOLA=10^{-9}$, $10^{-6}$ and $10^{-3}$,
\item GMRES with algebraic goal-oriented stopping criterion \eqref{aDWR}
  with tolerances $\cA=10^{-2}$, $10^{-1}$ and $10^{0}$,
\item BiCG with algebraic goal-oriented stopping criterion \eqref{aDWR}
  with tolerances $\cA=10^{-2}$, $10^{-1}$ and $10^{0}$,
\item BiCG with $\sigma$-stopping criterion \eqref{xiB}
  with tolerances $\cA=10^{-2}$, $10^{-1}$ and $10^{0}$.
\end{enumerate}
For all solvers we use the block-ILU(0) preconditioner which is suitable
for discontinuous Galerkin method, cf. \cite{lin_algeb}.
For GMRES, we solve first the primal problem and then the dual problem,
for BiCG, both problems are solved at once, of course.
Moreover, GMRES is restarted always after 45 iterations.
In case (c), the  stopping criterion \eqref{aDWR} is tested
after one restart of GMRES for the primal as well as the dual problems.
In case (d),  the  stopping criterion \eqref{aDWR} is tested after 100
BiCG iterations.

The results are presented in Tables~\ref{tab:CRh} -- \ref{tab:JBhp},
where we show the number of degrees of freedom $\DoF$ on the last mesh,
the total error $\J(\u-\uhmk)$, the final estimator of the discretization and algebraic errors
$\etaSmk$ and $\etaAmk$, respectively and their sum $\etamk$.
Moreover, these tables contains the total number of (GMRES or BiCG) iterations
on all mesh levels $m=0,1,\dots$ (iters)
and the total computational time in seconds.
If $\etaAmk \ge \etaSmk$ then the line ends with the character '!'.

We remind that the error estimator $\etamk$ does not give the upper bound,
see Remark~\ref{rem:guaranteed}.
There are small differences for each of four examples but we can state the following
observations. 
\begin{itemize}
\item the use of the preconditioned residual stopping criterion \eqref{RESP} for
  GMRES as well as BiCG with too large $\TOLA$ leads to the dominance of the
  algebraic error and consequently $J(u-\uhmk)$ and $\etamk$ are much larger
  then the tolerance $\TOL$ from \eqref{goal} and \eqref{STOP}, cf. (a) and (b),
\item GMRES requires more iterations  then BiCG
  since the primal and dual problems are solver separately, therefore
  BiCG is faster, compare (a) vs. (b) and (c) vs. (d),
\item  the goal-oriented stopping criterion \eqref{aDWR} allows to control
  the algebraic (and therefore also the total) error as expected for GMRES and BiCG,
  see (c) and (d),
\item the $\sigma$-stopping criterion \eqref{xiB} works efficiently, the algebraic error
  estimator $\etaAmk$ is smaller than the discretization one $\etaSmk$,
  the computational time is substantially reduced in comparison with
  the goal-oriented criterion \eqref{aDWR} since the higher-order reconstructions
  is performed less frequently, see the last paragraph in Section~\ref{sec:stop1}.
\end{itemize}

\begin{table}
  \bgroup
  \def\arraystretch{1.1}
   \setlength{\tabcolsep}{4pt}
 \begin{tabular}{c|ccccc|rrl}
 \hline
  $\TOLA$ or $\cA$  & $\DoF$ & $\J(\u-\uhmk)$ & $\etaSmk$ & $\etaAmk$ & $\etamk$  & iters & time(s) \\
 \hline
 \hline
   \multicolumn{5}{l}{(a) \ \      GMRES      preconditioned residual stopping criterion                  }  \\
 \hline
    1.E-09& 92680&   2.33E-08&   6.72E-09&   8.80E-11&   6.81E-09&    8301&      480.5&\!\!\!        \\ 
    1.E-06& 93120&   9.76E-08&   7.47E-09&   5.24E-08&   5.99E-08&    3274&      235.0&\!\!\!{\bf !} \\ 
    1.E-03& 94680&   1.27E-04&   7.36E-09&   9.25E-05&   9.25E-05&     806&      144.2&\!\!\!{\bf !} \\ 
 \hline
   \multicolumn{5}{l}{(b) \ \      BiCG       preconditioned residual stopping criterion                  }  \\
 \hline
    1.E-09& 92680&   2.32E-08&   6.72E-09&   6.26E-13&   6.72E-09&    3150&      341.8&\!\!\!        \\ 
    1.E-06& 90700&   2.67E-08&   7.96E-09&   5.07E-11&   8.01E-09&    2170&      258.6&\!\!\!        \\ 
    1.E-03& 90050&   4.80E-07&   7.94E-09&   2.54E-07&   2.62E-07&     900&      160.0&\!\!\!{\bf !} \\ 
 \hline
   \multicolumn{5}{l}{(c) \ \      GMRES      algebraic goal-oriented stopping criterion                  }  \\
 \hline
    1.E-02& 91770&   3.19E-08&   9.46E-09&   3.05E-11&   9.49E-09&    7560&      867.2&\!\!\!        \\ 
    1.E-01& 95230&   2.61E-08&   6.74E-09&   1.77E-10&   6.92E-09&    6480&      768.5&\!\!\!        \\ 
    1.E+00& 93440&   3.06E-08&   7.09E-09&   2.12E-09&   9.21E-09&    4860&      553.7&\!\!\!        \\ 
 \hline
   \multicolumn{5}{l}{(d) \ \      BiCG       algebraic goal-oriented stopping criterion                  }  \\
 \hline
    1.E-02& 94170&   2.56E-08&   6.71E-09&   1.19E-11&   6.72E-09&    3400&      677.9&\!\!\!        \\ 
    1.E-01& 90800&   2.69E-08&   7.83E-09&   5.76E-10&   8.40E-09&    2700&      534.8&\!\!\!        \\ 
    1.E+00& 91680&   3.67E-08&   9.34E-09&   2.29E-09&   1.16E-08&    2200&      374.3&\!\!\!        \\ 
 \hline
   \multicolumn{5}{l}{(e) \ \      BiCG       $\zeta$-stopping criterion                                  }  \\
 \hline
    1.E-02& 89760&   8.20E-06&   8.42E-09&   6.13E-06&   6.13E-06&     750&      233.3&\!\!\!{\bf !} \\ 
    1.E-01& 94320&   1.74E-05&   6.07E-09&   1.16E-05&   1.16E-05&     610&      179.5&\!\!\!{\bf !} \\ 
    1.E+00& 91660&   5.77E-05&   8.22E-09&   4.35E-05&   4.35E-05&     430&      160.1&\!\!\!{\bf !} \\ 
 \hline
   \multicolumn{5}{l}{(f) \ \      BiCG       $\eta$-stopping criterion                                   }  \\
 \hline
    1.E-02&122070&   1.05E-08&   3.14E-09&   2.84E-11&   3.17E-09&    2950&      479.7&\!\!\!        \\ 
    1.E-01&121760&   8.78E-09&   2.45E-09&   7.33E-10&   3.19E-09&    2550&      399.7&\!\!\!        \\ 
    1.E+00&121620&   1.42E-08&   3.16E-09&   1.82E-09&   4.99E-09&    2270&      411.3&\!\!\!        \\ 
 \hline
   \multicolumn{5}{l}{(g) \ \      BiCG       $\sigma$-stopping criterion                                 }  \\
 \hline
    1.E-02& 93560&   2.59E-08&   6.90E-09&   2.41E-11&   6.93E-09&    2340&      277.4&\!\!\!        \\ 
    1.E-01& 90730&   2.94E-08&   7.68E-09&   2.77E-10&   7.96E-09&    1890&      229.7&\!\!\!        \\ 
    1.E+00& 91110&   3.40E-08&   7.30E-09&   5.15E-09&   1.25E-08&    1730&      226.5&\!\!\!        \\ 
 \hline
 \end{tabular}

  \egroup
  \caption{Elliptic problem \eqref{eq:CR}, $h$-mesh adaptation with $P_3$ approximation,
    comparison of solvers and stopping criteria, $\TOL=10^{-8}$.}
  \label{tab:CRh}  
\end{table}

\begin{table}
  \bgroup
  \def\arraystretch{1.1}
   \setlength{\tabcolsep}{4pt}
 \begin{tabular}{c|ccccc|rrl}
 \hline
  $\TOLA$ or $\cA$  & $\DoF$ & $\J(\u-\uhmk)$ & $\etaSmk$ & $\etaAmk$ & $\etamk$  & iters & time(s) \\
 \hline
 \hline
   \multicolumn{5}{l}{(a) \ \      GMRES      preconditioned residual stopping criterion                  }  \\
 \hline
    1.E-09& 22390&   3.15E-08&   6.02E-09&   9.33E-12&   6.03E-09&    4810&      183.1&\!\!\!        \\ 
    1.E-06& 24160&   3.84E-08&   5.17E-09&   2.66E-08&   3.17E-08&    2627&      121.4&\!\!\!{\bf !} \\ 
    1.E-03& 23845&   2.90E-05&   8.24E-09&   1.83E-05&   1.83E-05&     602&       87.8&\!\!\!{\bf !} \\ 
 \hline
   \multicolumn{5}{l}{(b) \ \      BiCG       preconditioned residual stopping criterion                  }  \\
 \hline
    1.E-09& 23520&   3.38E-08&   7.73E-09&   7.46E-14&   7.73E-09&    2360&      152.0&\!\!\!        \\ 
    1.E-06& 25656&   2.07E-08&   3.80E-09&   6.73E-11&   3.87E-09&    1700&      138.4&\!\!\!        \\ 
    1.E-03& 25364&   9.51E-08&   7.03E-09&   3.15E-08&   3.85E-08&     970&      114.5&\!\!\!{\bf !} \\ 
 \hline
   \multicolumn{5}{l}{(c) \ \      GMRES      algebraic goal-oriented stopping criterion                  }  \\
 \hline
    1.E-02& 24575&   2.58E-08&   5.00E-09&   6.44E-13&   5.00E-09&    3753&      200.7&\!\!\!        \\ 
    1.E-01& 24468&   2.70E-08&   5.40E-09&   1.09E-10&   5.51E-09&    3408&      167.3&\!\!\!        \\ 
    1.E+00& 23594&   2.75E-08&   5.22E-09&   1.13E-10&   5.34E-09&    3228&      158.5&\!\!\!        \\ 
 \hline
   \multicolumn{5}{l}{(d) \ \      BiCG       algebraic goal-oriented stopping criterion                  }  \\
 \hline
    1.E-02& 23631&   3.12E-08&   6.41E-09&   3.96E-12&   6.41E-09&    2100&      186.8&\!\!\!        \\ 
    1.E-01& 23721&   2.89E-08&   5.14E-09&   2.57E-11&   5.17E-09&    1900&      163.6&\!\!\!        \\ 
    1.E+00& 24006&   2.54E-08&   4.80E-09&   1.23E-12&   4.80E-09&    1800&      162.1&\!\!\!        \\ 
 \hline
   \multicolumn{5}{l}{(e) \ \      BiCG       $\zeta$-stopping criterion                                  }  \\
 \hline
    1.E-02& 26627&   1.07E-07&   5.36E-09&   6.34E-08&   6.87E-08&     830&      104.7&\!\!\!{\bf !} \\ 
    1.E-01& 22219&   2.33E-06&   9.81E-09&   1.24E-06&   1.25E-06&     600&       89.5&\!\!\!{\bf !} \\ 
    1.E+00& 31772&   1.12E-05&   8.16E-09&   1.24E-05&   1.24E-05&     520&      139.8&\!\!\!{\bf !} \\ 
 \hline
   \multicolumn{5}{l}{(f) \ \      BiCG       $\eta$-stopping criterion                                   }  \\
 \hline
    1.E-02& 24783&   1.93E-08&   3.94E-09&   2.42E-11&   3.96E-09&    1820&      168.3&\!\!\!        \\ 
    1.E-01& 26703&   1.41E-08&   2.25E-09&   1.07E-10&   2.35E-09&    1680&      179.3&\!\!\!        \\ 
    1.E+00& 26491&   2.85E-08&   4.84E-09&   1.83E-09&   6.67E-09&    1380&      165.6&\!\!\!        \\ 
 \hline
   \multicolumn{5}{l}{(g) \ \      BiCG       $\sigma$-stopping criterion                                 }  \\
 \hline
    1.E-02& 24176&   2.49E-08&   5.07E-09&   3.08E-11&   5.10E-09&    1660&      128.0&\!\!\!        \\ 
    1.E-01& 23064&   2.84E-08&   5.69E-09&   6.39E-10&   6.33E-09&    1460&      117.4&\!\!\!        \\ 
    1.E+00& 22486&   2.84E-08&   6.04E-09&   9.92E-10&   7.04E-09&    1300&      112.9&\!\!\!        \\ 
 \hline
 \end{tabular}

  \egroup
  \caption{Elliptic problem \eqref{eq:CR}, $hp$-mesh adaptation,
    comparison of solvers and stopping criteria, $\TOL=10^{-8}$.}
  \label{tab:CRhp}  
\end{table}

\begin{table}
  \bgroup
  \def\arraystretch{1.1}
   \setlength{\tabcolsep}{4pt}
 \begin{tabular}{c|ccccc|rrl}
 \hline
  $\TOLA$ or $\cA$  & $\DoF$ & $\J(\u-\uhmk)$ & $\etaSmk$ & $\etaAmk$ & $\etamk$  & iters & time(s) \\
 \hline
 \hline
   \multicolumn{5}{l}{(a) \ \      GMRES      preconditioned residual stopping criterion                  }  \\
 \hline
    1.E-09& 74910&   4.36E-12&   6.09E-11&   2.55E-10&   3.15E-10&    4862&      308.2&\!\!\!{\bf !} \\ 
    1.E-06& 83960&   7.36E-08&   3.22E-11&   5.29E-08&   5.29E-08&    2198&      240.5&\!\!\!{\bf !} \\ 
    1.E-03& 70060&   4.16E-06&   3.25E-11&   1.35E-04&   1.35E-04&     524&      163.5&\!\!\!{\bf !} \\ 
 \hline
   \multicolumn{5}{l}{(b) \ \      BiCG       preconditioned residual stopping criterion                  }  \\
 \hline
    1.E-09& 70580&   1.07E-10&   6.37E-11&   6.79E-12&   7.05E-11&    2920&      327.1&\!\!\!        \\ 
    1.E-06& 78610&   4.95E-11&   2.80E-11&   4.39E-10&   4.67E-10&    2310&      317.1&\!\!\!{\bf !} \\ 
    1.E-03& 61250&   1.37E-08&   1.11E-12&   2.01E-07&   2.01E-07&    1100&      195.5&\!\!\!{\bf !} \\ 
 \hline
   \multicolumn{5}{l}{(c) \ \      GMRES      algebraic goal-oriented stopping criterion                  }  \\
 \hline
    1.E-02& 52820&   3.04E-10&   8.52E-11&   2.49E-13&   8.54E-11&    5501&      430.6&\!\!\!        \\ 
    1.E-01& 79990&   6.36E-11&   5.80E-11&   2.99E-12&   6.10E-11&    6404&      679.7&\!\!\!        \\ 
    1.E+00& 68850&   1.50E-10&   6.06E-11&   1.49E-11&   7.56E-11&    5324&      561.3&\!\!\!        \\ 
 \hline
   \multicolumn{5}{l}{(d) \ \      BiCG       algebraic goal-oriented stopping criterion                  }  \\
 \hline
    1.E-02& 52680&   2.46E-10&   8.06E-11&   4.13E-14&   8.07E-11&    4180&      601.8&\!\!\!        \\ 
    1.E-01& 73320&   1.15E-10&   4.66E-11&   3.65E-14&   4.66E-11&    4180&      694.3&\!\!\!        \\ 
    1.E+00& 79170&   8.30E-11&   3.09E-11&   9.52E-12&   4.04E-11&    4780&     1074.4&\!\!\!        \\ 
 \hline
   \multicolumn{5}{l}{(e) \ \      BiCG       $\zeta$-stopping criterion                                  }  \\
 \hline
    1.E-02& 76330&   3.32E-07&   4.21E-11&   9.32E-04&   9.32E-04&     880&      208.1&\!\!\!{\bf !} \\ 
    1.E-01&106050&   4.05E-06&   1.41E-11&   1.87E-03&   1.87E-03&     680&      214.5&\!\!\!{\bf !} \\ 
    1.E+00& 63940&   3.12E-06&   9.97E-11&   1.37E-03&   1.37E-03&     560&      152.1&\!\!\!{\bf !} \\ 
 \hline
   \multicolumn{5}{l}{(f) \ \      BiCG       $\eta$-stopping criterion                                   }  \\
 \hline
    1.E-02& 72820&   8.49E-11&   2.74E-11&   4.50E-13&   2.78E-11&    3070&      398.8&\!\!\!        \\ 
    1.E-01& 71380&   9.39E-11&   3.25E-11&   2.54E-12&   3.50E-11&    2840&      433.1&\!\!\!        \\ 
    1.E+00& 81880&   5.95E-11&   2.53E-11&   1.28E-11&   3.81E-11&    2630&      468.1&\!\!\!        \\ 
 \hline
   \multicolumn{5}{l}{(g) \ \      BiCG       $\sigma$-stopping criterion                                 }  \\
 \hline
    1.E-02& 72820&   8.49E-11&   2.74E-11&   4.50E-13&   2.78E-11&    3070&      357.6&\!\!\!        \\ 
    1.E-01& 71380&   9.39E-11&   3.25E-11&   2.54E-12&   3.50E-11&    2840&      319.4&\!\!\!        \\ 
    1.E+00& 81880&   5.95E-11&   2.53E-11&   1.28E-11&   3.81E-11&    2630&      338.7&\!\!\!        \\ 
 \hline
 \end{tabular}

  \egroup
  \caption{Convection-diffusion problem \eqref{carpio},  $h$-mesh adaptation with
    $P_3$-approximation,  comparison of solvers and stopping criteria, $\TOL=10^{-10}$.}
  \label{tab:JBP3}  
\end{table}

\begin{table}
  \bgroup
  \def\arraystretch{1.1}
   \setlength{\tabcolsep}{4pt}
 \begin{tabular}{c|ccccc|rrl}
 \hline
  $\TOLA$ or $\cA$  & $\DoF$ & $\J(\u-\uhmk)$ & $\etaSmk$ & $\etaAmk$ & $\etamk$  & iters & time(s) \\
 \hline
 \hline
   \multicolumn{5}{l}{(a) \ \      GMRES      preconditioned residual stopping criterion                  }  \\
 \hline
    1.E-09& 15852&   5.93E-10&   4.54E-13&   6.76E-11&   6.81E-11&    2909&      124.0&\!\!\!{\bf !} \\ 
    1.E-06& 18233&   1.46E-09&   6.76E-11&   3.33E-08&   3.33E-08&    2006&      151.9&\!\!\!{\bf !} \\ 
    1.E-03& 20422&   3.08E-06&   4.89E-11&   9.65E-06&   9.65E-06&     590&      119.0&\!\!\!{\bf !} \\ 
 \hline
   \multicolumn{5}{l}{(b) \ \      BiCG       preconditioned residual stopping criterion                  }  \\
 \hline
    1.E-09& 19905&   6.20E-11&   4.33E-11&   7.96E-12&   5.13E-11&    1980&      178.0&\!\!\!        \\ 
    1.E-06& 18408&   2.26E-11&   5.13E-11&   1.39E-10&   1.91E-10&    1550&      163.6&\!\!\!{\bf !} \\ 
    1.E-03& 23431&   2.56E-08&   2.50E-11&   1.49E-06&   1.49E-06&     820&      146.0&\!\!\!{\bf !} \\ 
 \hline
   \multicolumn{5}{l}{(c) \ \      GMRES      algebraic goal-oriented stopping criterion                  }  \\
 \hline
    1.E-02& 19323&   5.46E-11&   3.29E-12&   5.32E-15&   3.29E-12&    4690&      263.9&\!\!\!        \\ 
    1.E-01& 17504&   1.81E-10&   9.99E-11&   1.66E-14&   9.99E-11&    3876&      221.3&\!\!\!        \\ 
    1.E+00& 24723&   1.25E-11&   4.92E-11&   4.30E-13&   4.97E-11&    4236&      282.6&\!\!\!        \\ 
 \hline
   \multicolumn{5}{l}{(d) \ \      BiCG       algebraic goal-oriented stopping criterion                  }  \\
 \hline
    1.E-02& 18689&   8.64E-11&   5.22E-12&   2.18E-14&   5.25E-12&    3870&      383.5&\!\!\!        \\ 
    1.E-01& 18668&   9.03E-11&   2.44E-11&   3.54E-13&   2.48E-11&    2670&      257.1&\!\!\!        \\ 
    1.E+00& 17240&   2.64E-10&   5.22E-11&   1.61E-12&   5.38E-11&    2270&      203.2&\!\!\!        \\ 
 \hline
   \multicolumn{5}{l}{(e) \ \      BiCG       $\zeta$-stopping criterion                                  }  \\
 \hline
    1.E-02& 20038&   3.09E-08&   1.35E-11&   2.25E-05&   2.25E-05&     880&      176.6&\!\!\!{\bf !} \\ 
    1.E-01& 21598&   9.45E-07&   5.14E-11&   1.99E-04&   1.99E-04&     780&      170.5&\!\!\!{\bf !} \\ 
    1.E+00& 34233&   2.72E-04&   6.20E-11&   1.40E-03&   1.40E-03&     810&      161.3&\!\!\!{\bf !} \\ 
 \hline
   \multicolumn{5}{l}{(f) \ \      BiCG       $\eta$-stopping criterion                                   }  \\
 \hline
    1.E-02& 20146&   2.92E-10&   2.55E-11&   2.89E-13&   2.58E-11&    1910&      181.2&\!\!\!        \\ 
    1.E-01& 22752&   7.32E-11&   3.73E-11&   1.06E-12&   3.83E-11&    2090&      209.8&\!\!\!        \\ 
    1.E+00& 22936&   4.23E-12&   2.00E-11&   2.32E-12&   2.23E-11&    1900&      228.2&\!\!\!        \\ 
 \hline
   \multicolumn{5}{l}{(g) \ \      BiCG       $\sigma$-stopping criterion                                 }  \\
 \hline
    1.E-02& 20146&   2.92E-10&   2.55E-11&   2.89E-13&   2.58E-11&    1910&      156.3&\!\!\!        \\ 
    1.E-01& 18757&   2.05E-10&   1.92E-11&   1.32E-12&   2.05E-11&    1790&      148.7&\!\!\!        \\ 
    1.E+00& 22936&   4.23E-12&   2.00E-11&   2.32E-12&   2.23E-11&    1900&      184.0&\!\!\!        \\ 
 \hline
 \end{tabular}

  \egroup
  \caption{Convection-diffusion problem \eqref{carpio}, $hp$-mesh adaptation,
    comparison of solvers and stopping criteria, $\TOL=10^{-10}$.}
  \label{tab:JBhp}  
\end{table}

\section{Summary of the results and outlook}
\label{sec:sum}

We developed an efficient technique for the numerical solution of
primal and dual algebraic systems arising in the goal-oriented
error estimation and mesh adaptation.
Both algebraic systems are solved simultaneously by BiCG method which allows
to control the algebraic error during the iterative process.
The proposed $\sigma$-stopping criterion is cheap for the evaluation and
significantly reduce the computational costs.
Moreover, it guarantees that the algebraic error estimate
bounded by the discretization one.

Further natural step is to extend this approach for the solution of nonlinear
partial differential equations. The the dual problem has to be build on a
linearization of the primal one. However, employing a Newton-like method
for the solution of the discretized primal problem,
an approximate solution of the dual problem is available
at each Newton step 
and the technique developed in this paper can be employed.
However, it is necessary to balance
the linear algebraic errors, 
the non-linear algebraic errors and
the discretization errors. This is the subject of the further work.

\paragraph{Acknowledgements} 
The authors are thankful to their colleagues from the Char\-les University,
namely M.Kub{\'\i}nov{\'a}, T. Gergelits and F. Roskovec
for a fruitful discussion.

\bibliographystyle{spmpsci}      
\bibliography{ref}   

\begin{thebibliography}{10}
\providecommand{\url}[1]{{#1}}
\providecommand{\urlprefix}{URL }
\expandafter\ifx\csname urlstyle\endcsname\relax
  \providecommand{\doi}[1]{DOI~\discretionary{}{}{}#1}\else
  \providecommand{\doi}{DOI~\discretionary{}{}{}\begingroup
  \urlstyle{rm}\Url}\fi

\bibitem{AinsworthOden}
Ainsworth, M., Oden, J.T.: A posteriori error estimation in finite element
  analysis.
\newblock Pure and Applied Mathematics (New York). Wiley-Interscience [John
  Wiley \& Sons], New York (2000)

\bibitem{Ainsworth2012Guaranteed}
Ainsworth, M., Rankin, R.: {Guaranteed computable bounds on quantities of
  interest in finite element computations}.
\newblock International Journal for Numerical Methods in Engineering
  \textbf{89}(13), 1605--1634 (2012)

\bibitem{Arioli_NM04}
Arioli, M.: A stopping criterion for the conjugate gradient algorithm in a
  finite element method framework.
\newblock Numer. Math. \textbf{97}(1), 1--24 (2004)

\bibitem{ArioliLiiesenMiedlerStrakos13}
Arioli, M., Liesen, J., Mi{e}dlar, A., Strako{\v{s}}, Z.: Interplay between
  discretization and algebraic computation in adaptive numerical solution of
  elliptic {PDE} problems.
\newblock GAMM-Mitt. \textbf{36}(1), 102--129 (2013)

\bibitem{BalanWoopenMay16}
Balan, A., Woopen, M., May, G.: Adjoint-based $hp$-adaptivity on anisotropic
  meshes for high-order compressible flow simulations.
\newblock Comput. Fluids \textbf{139}, 47 -- 67 (2016)

\bibitem{RannacherBook}
Bangerth, W., Rannacher, R.: {Adaptive Finite Element Methods for Differential
  Equations}.
\newblock {Lectures in Mathematics. ETH Z{\"u}rich}. {Birkh\"auser Verlag}
  (2003)

\bibitem{B:BaBeCh1994}
Barrett, R., Berry, M., Chan, T.F., et~al.: Templates for the solution of
  linear systems: building blocks for iterative methods.
\newblock Society for Industrial and Applied Mathematics (SIAM), Philadelphia,
  PA (1994).
\newblock \doi{10.1137/1.9781611971538}.
\newblock \urlprefix\url{https://doi.org/10.1137/1.9781611971538}

\bibitem{ESCO-18}
Barto{\v s}, O., Dolej{\v s}{\'\i}, V., May, G., , Rangarajan, A., Roskovec,
  F.: Goal-oriented anisotropic $hp$-mesh optimization technique for linear
  convection-diffusion-reaction problem.
\newblock Comput. Math. Appl. \textbf{78}(9), 2973--2993 (2019)

\bibitem{BeckerRannacher01}
Becker, R., Rannacher, R.: An optimal control approach to a-posteriori error
  estimation in finite element methods.
\newblock Acta Numerica \textbf{10}, 1--102 (2001)

\bibitem{CarpioPrietoBermejo_SISC13}
Carpio, J., Prieto, J., Bermejo, R.: Anisotropic "goal-oriented" mesh
  adaptivity for elliptic problems.
\newblock SIAM J. Sci. Comput. \textbf{35}(2), A861--A885 (2013)

\bibitem{ADGFEM}
Dolej{\v s}{\'\i}, V.: ADGFEM -- Adaptive discontinuous {G}alerkin finite
  element method, in-house code.
\newblock Charles University, Prague, Faculty of Mathematics and Physics
  (2014).
\newblock \url{http://atrey.karlin.mff.cuni.cz/~dolejsi/adgfem/}

\bibitem{DGM-book}
Dolej{\v s}{\'\i}, V., Feistauer, M.: Discontinuous Galerkin Method -- Analysis
  and Applications to Compressible Flow.
\newblock Springer Series in Computational Mathematics 48. Springer, Cham
  (2015)

\bibitem{lin_algeb}
Dolej{\v s}{\'\i}, V., Hol{\'\i}k, M., Hozman, J.: Efficient solution strategy
  for the semi-implicit discontinuous {G}alerkin discretization of the
  {N}avier-{S}tokes equations.
\newblock J. Comput. Phys. \textbf{230}, 4176--4200 (2011)

\bibitem{DolRoskovec_AM17}
Dolej{\v s}{\'\i}, V., Roskovec, F.: Goal-oriented error estimates including
  algebraic errors in discontinuous {G}alerkin discretizations of linear
  boundary value problems.
\newblock Appl. Math. \textbf{62}(6), 579--605 (2017)

\bibitem{FidkowskiDarmofal_AIAA11}
Fidkowski, K., Darmofal, D.: Review of output-based error estimation and mesh
  adaptation in computational fluid dynamics.
\newblock AIAA Journal \textbf{49}(4), 673--694 (2011)

\bibitem{Fl1976}
Fletcher, R.: Conjugate gradient methods for indefinite systems.
\newblock In: Numerical analysis (Proc 6th Biennial Dundee Conf., Univ. Dundee,
  Dundee, 1975), pp. 73--89. Lecture Notes in Math., Vol. 506. Springer, Berlin
  (1976).
\newblock \urlprefix\url{http://www.ams.org/mathscinet-getitem?mr=57\#1841}

\bibitem{FormaggiaPerotto04}
Formaggia, L., Micheletti, S., Perotto, S.: Anisotropic mesh adaption in
  computational fluid dynamics: application to the advection-diffusion-reaction
  and the {S}tokes problems.
\newblock Appl. Numer. Math. \textbf{51}(4), 511--533 (2004)

\bibitem{GeHaHo09}
Georgoulis, E.H., Hall, E., Houston, P.: Discontinuous {G}alerkin methods on
  $hp$-anisotropic meshes {II}: {A} posteriori error analysis and adaptivity.
\newblock Appl. Numer. Math. \textbf{59}(9), 2179--2194 (2009)

\bibitem{GileSuli02}
Giles, M., S{\"u}li, E.: Adjoint methods for {PDE}s: a posteriori error
  analysis and postprocessing by duality.
\newblock Acta Numerica \textbf{11}, 145--236 (2002)

\bibitem{Gr1997}
Greenbaum, A.: Estimating the attainable accuracy of recursively computed
  residual methods.
\newblock SIAM J. Matrix Anal. Appl. \textbf{18}(3), 535--551 (1997)

\bibitem{Hartman08}
Hartmann, R.: Multitarget error estimation and adaptivity in aerodynamic flow
  simulations.
\newblock SIAM J. Sci. Comput. \textbf{31}(1), 708--731 (2008)

\bibitem{HH06:SIPG2}
Hartmann, R., Houston, P.: Symmetric interior penalty {DG} methods for the
  compressible {N}avier-{S}tokes equations {II}: {G}oal-oriented a posteriori
  error estimation.
\newblock Int. J. Numer. Anal. Model. \textbf{3}, 141--162 (2006)

\bibitem{VohralikStrakos10}
Jir\'anek, P., Strako{\v s}, Z., Vohral{\'\i}k, M.: A posteriori error
  estimates including algebraic error and stopping criteria for iterative
  solvers.
\newblock SIAM J. Sci. Comput. \textbf{32}(3), 1567--1590 (2010)

\bibitem{Korotov_JCAM06}
Korotov, S.: A posteriori error estimation of goal-oriented quantities for
  elliptic type {BVP}s.
\newblock J. Comput. Appl. Math. \textbf{191}(2), 216--227 (2006)

\bibitem{KuzminKorotov_MCS10}
Kuzmin, D., Korotov, S.: {Goal-oriented a posteriori error estimates for
  transport problems}.
\newblock {Math. Comput. Simul.} \textbf{{80}}({8, SI}), {1674--1683} ({2010})

\bibitem{KuzminMoller_JCAM10}
Kuzmin, D., M{\"o}ller, M.: Goal-oriented mesh adaptation for flux-limited
  approximations to steady hyperbolic problems.
\newblock J. Comput. Appl. Math. \textbf{233}(12), 3113--3120 (2010)

\bibitem{LoseilleDervieuxAlauzet_JCP10}
Loseille, A., Dervieux, A., Alauzet, F.: Fully anisotropic goal-oriented mesh
  adaptation for 3{D} steady {E}uler equations.
\newblock J. Comput. Phys. \textbf{229}(8), 2866--2897 (2010)

\bibitem{MallikVohralikYousef_JCAM20}
Mallik, G., Vohral{\'\i}k, M., Yousef, S.: Goal-oriented a posteriori error
  estimation for conforming and nonconforming approximations with inexact
  solvers.
\newblock Journal of Computational and Applied Mathematics \textbf{366} (2020)

\bibitem{Meidner2009Goal}
Meidner, D., Rannacher, R., Vihharev, J.: {Goal-oriented error control of the
  iterative solution of finite element equations}.
\newblock J. Numer. Math. \textbf{17}, 143 (2009)

\bibitem{NochettoVeeserVerani_IMAJNA09}
Nochetto, R., Veeser, A., Verani, M.: A safeguarded dual weighted residual
  method.
\newblock IMA Journal of Numerical Analysis \textbf{29}(1), 126--140 (2009)

\bibitem{OdenPrudhomme_CAMWA01}
Oden, J., Prudhomme, S.: Goal-oriented error estimation and adaptivity for the
  finite element method.
\newblock Comput. Math. Appl. \textbf{41}(5-6), 735--756 (2001)

\bibitem{Picasso_CNME09}
Picasso, M.: A stopping criterion for the conjugate gradient algorithm in the
  framework of anisotropic adaptive finite elements.
\newblock Communications in Numerical Methods in Engineering \textbf{25}(4),
  339--355 (2009)

\bibitem{Rey2Gosselet_CMAME15}
Rey, V., Gosselet, P., Rey, C.: Strict bounding of quantities of interest in
  computations based on domain decomposition.
\newblock Computer Methods in Applied Mechanics and Engineering \textbf{287}
  (2015)

\bibitem{Rey2Gosselet_CMAME14}
Rey, V., Rey, C., Gosselet, P.: A strict error bound with separated
  contributions of the discretization and of the iterative solver in
  non-overlapping domain decomposition methods.
\newblock Computer Methods in Applied Mechanics and Engineering \textbf{270},
  293--303 (2014)

\bibitem{Richter_IJNMF10}
Richter, T.: A posteriori error estimation and anisotropy detection with the
  dual-weighted residual method.
\newblock Int. J. Numer. Meth. Fluids \textbf{62}, 90--118 (2010)

\bibitem{SolDemko04}
{\v S}ol{\'\i}n, P., Demkowicz, L.: Goal-oriented $hp$-adaptivity for elliptic
  problems.
\newblock Comput. Methods Appl. Mech. Engrg. \textbf{193}, 449--468 (2004)

\bibitem{StTi2011}
Strako{\v s}, Z., Tich{\'y}, P.: On efficient numerical approximation of the
  bilinear form $c^{\star}{A}^{-1}b$.
\newblock SIAM J. Sci. Comput. \textbf{33}(2), 565--587 (2010)

\bibitem{verfurth-book2}
Verf\"urth, R.: A Posteriori Error Estimation Techniques for Finite Element
  Methods.
\newblock Numerical Mathematics and Scientific Computation. Oxford University
  Press (2013)

\end{thebibliography}

\end{document}